\documentclass[10pt]{article}
\usepackage{amsmath,amssymb}
\usepackage[dvips]{graphicx}
\usepackage{theorem}
\usepackage{esint}
\usepackage{color}
\usepackage{booktabs}
\usepackage{multirow}
\usepackage{bm}
\usepackage{soul}

\usepackage{epstopdf}
\usepackage{mathtools}

\numberwithin{equation}{section}{\theorembodyfont{\upshape}
\newtheorem{remark}{\textbf{\emph{Remark}}}[section]

\theoremstyle{plain}
\theoremstyle{plain} 
\theoremstyle{plain}
\numberwithin{equation}{section}
{\theorembodyfont{\upshape}\theoremstyle{plain}}

\newcommand{\ih}{\mathcal{I}_h}
\newcommand{\ov}{\overline}
\newcommand{\mc}{\mathcal}

\newcommand{\bi}{\bm{i}}

\usepackage[left=3.5cm,right=3cm]{geometry}
\usepackage{lineno}

\newcommand*\patchAmsMathEnvironmentForLineno[1]{%
  \expandafter\let\csname old#1\expandafter\endcsname\csname #1\endcsname
  \expandafter\let\csname oldend#1\expandafter\endcsname\csname end#1\endcsname
  \renewenvironment{#1}%
     {\linenomath\csname old#1\endcsname}%
     {\csname oldend#1\endcsname\endlinenomath}}%
\newcommand*\patchBothAmsMathEnvironmentsForLineno[1]{%
  \patchAmsMathEnvironmentForLineno{#1}%
  \patchAmsMathEnvironmentForLineno{#1*}}%
\AtBeginDocument{%
\patchBothAmsMathEnvironmentsForLineno{equation}%
\patchBothAmsMathEnvironmentsForLineno{align}%
\patchBothAmsMathEnvironmentsForLineno{flalign}%
\patchBothAmsMathEnvironmentsForLineno{alignat}%
\patchBothAmsMathEnvironmentsForLineno{gather}%
\patchBothAmsMathEnvironmentsForLineno{multline}%
}

\usepackage[]{algorithm2e}

\graphicspath{{figures/}}

    \makeatletter
    \let\@fnsymbol\@arabic
    \makeatother

\begin{document}

\title{Stable Generalized Finite Element Method and
associated iterative schemes; application to interface problems}
\author{Kenan Kergrene
\thanks{
Department of Mathematics and Industrial Engineering, \'{E}cole Polytechnique de Montr\'{e}al, Canada.
} 
\and
Ivo Babu\v{s}ka \thanks{
ICES, University of Texas at Austin, Austin, TX, United States.
}
\and Uday Banerjee
\thanks{
Department of Mathematics, 215 Carnegie, Syracuse University,
Syracuse, NY 13244, United States. E-mail: banerjee@syr.edu.} 
}

\date{}

\maketitle

\begin{abstract}
The Generalized Finite Element Method (GFEM) is an extension of the Finite Element Method (FEM), where the standard finite element space is augmented with a space of non-polynomial functions, called the enrichment space. The functions in the enrichment space mimic the local behavior of the unknown solution of the underlying variational problem. GFEM has been successfully applied to a wide range of problems. However, it often suffers from bad conditioning, i.e., its conditioning may not be robust with respect to the mesh and in fact, the conditioning could be much worse than that of the standard FEM. In this paper, we present a numerical study that shows that if the ``angle" between the finite element space and the enrichment space is bounded away from 0, uniformly with respect to the mesh, then the GFEM is stable, i.e., the conditioning of GFEM is not worse than that of the standard FEM. A GFEM with this property is called a Stable GFEM (SGFEM). The last part of the paper is devoted to the derivation of a robust iterative solver exploiting this angle condition. It is shown that the required ``wall-clock" time is greatly reduced compared to popular GFEMs used in the literature.
\end{abstract}

\noindent {\bf Keywords: Generalized Finite Element Method (GFEM), Partition of Unity Method (PUM), Stable GFEM (SGFEM), Condition Number, Angle Condition} 

\setlength{\baselineskip}{15pt}

\section{Introduction}

The Generalized Finite Element Method (GFEM) has sparked a lot of interest in the last 20 years and has been successfully applied to a wide range of engineering problems, e.g., crack-propagation, material modeling, and solid--fluid interactions. We refer to the review articles \cite{AbazizHam,BelGraVen,EfenHou,FriBel} and the citations therein for various applications of GFEM. The method has been incorporated into commercial codes, e.g., Abaqus and LS-DYNA \cite{Abaqus,LSDYNA}. It is also known in the literature as the Extended Finite Element Method (XFEM). We will simply refer to the method as GFEM and we will address special instances of this method such as SGFEM and M-GFEM.

As the name GFEM/XFEM suggests, the GFEM is a generalization/extension of the standard Finite Element Method (FEM). Specific non-polynomial local basis functions that mimic special features (e.g., singularity) of the solution of the underlying PDE model of interest, are used in this method in addition to the standard ``hat-functions.'' These additional local basis functions are called the local \textit{enrichment} functions. In fact, the GFEM is a particular instance of the Partition of Unity Method (PUM). The PUM, developed in \cite{BabMel,Mel,MelBab}, allows the use of any Partition of Unity (PU) together with local enrichment functions. The GFEM is a PUM, where the finite element ``hat-functions'' serve as the PU. Various methods for solving multi-scale problems are also based on PUM; see for example \cite{BabHuaLip,BabLip,EfenHou}. The PUM based on a ``flat-top" PU was developed in \cite{GrSch7,Sch}. A similar idea, referred to as $h-p$ Cloud method, was developed in \cite{DuarOden,DuarOden2}. The original idea of GFEM, i.e., the use of hat-functions as the PU, was introduced in \cite{BabCalOsb}. Since then, the GFEM has been developed, refined, and used in various applications in two and three dimensions, e.g., in \cite{BelyBla,DolbowPhd,DuaBabOden2000,DuaHamLis,DuaOden1997,StrBabCop,StrCopBab1,StrCopBab,SukChopEtAl}. The GFEM is often referred to as the XFEM in the literature. We mention that the use of local non-polynomial approximation, not in the framework of PUM, was suggested earlier in \cite{BabOsb2}.

The XFEM/GFEM was initially developed as a computational method with essentially intuitive understanding of the necessity of appropriate enrichments for convergence. The method was primarily tested numerically to ensure convergence. Appropriate enrichment functions for various applications were identified in the literature and the optimal convergence of the approximate solution was shown through computations. A rigorous mathematical proof of optimal convergence was derived in \cite{NicRenard}, in the context of a crack problem.

Though approximability and optimal convergence are very important features of a numerical method such as GFEM, it is equally important that the underlying linear system could be solved accurately and efficiently. Solving such linear systems accurately and efficiently depends on the stability of the GFEM, i.e., on the conditioning of the underlying linear system. It was reported early in \cite{BabMel,FriBel} that the GFEM could be unstable and that its conditioning may not be robust with respect to the mesh. However, there are very few papers that addressed these issues by carefully studying the conditioning of the GFEM and by examining the performance of associated iterative solvers to solve the linear system. Various ad-hoc stabilization procedures were used to address these issues in \cite{BecMinEtAl,Loe,MenBor,StrBabCop,StrCopBab}. Stabilization based on a local orthogonalization idea was used in \cite{Loe}. We mention however that local orthogonalization was also addressed in the context of PUM with flat-top PU in \cite{Schw3}.

Extensive literature is available on the loss of accuracy in the computed solution of a linear system; we refer to the monographs \cite{Higham,Wilkinson}. In \cite{BabBan}, it was shown that the Scaled Condition Number (SCN) of the matrix, related to FEM and GFEM, is a good indicator of the stability and loss of accuracy in the solution obtained from elimination methods. 

The conditioning of GFEM was addressed in \cite{BabBan} where the idea of a stable GFEM (SGFEM) was introduced. In general, a GFEM is called stable (SGFEM) if 
\begin{itemize}
\item[(i)] it yields the optimal order of convergence, and
\item[(ii)] the SCN of the linear system associated with the GFEM is of the same order $O(h^{-2})$ ($h$ being the discretization parameter) as that of a standard FEM in a robust manner with respect to the mesh.
\end{itemize}
It was shown in \cite{BabBan} that the SCN of a GFEM could be much higher than that of the FEM, e.g., $O(h^{-4})$. It was mathematically established in that paper that if the enrichments satisfy two specific conditions, then the SCN of the underlying GFEM is of the same order as that of a standard FEM. For various problems in the 1-D setting, a simple modification based on subtracting the piecewise linear interpolant of the standard enrichment was suggested in \cite{BabBan} and it was shown that the modified GFEM was indeed an SGFEM for these problems. However, the modification suggested in \cite{BabBan} may lead to loss of accuracy in some problems in higher dimensions as shown in \cite{GupDuaEtAl,GupDua2EtAl,SauFri2}. It was shown in \cite{GupDuaEtAl,GupDua2EtAl} that a further modification of ``Heaviside enrichment," in the context of a problem with a crack, is required for a GFEM to be an SGFEM, i.e., the further modification restores the accuracy of the computed solution while retaining the well-conditioning of the linear system. Thus a GFEM with the simple modification of enrichments as suggested in \cite{BabBan} may not yield an SGFEM for every problem; further modifications of the enrichments may be required for a GFEM to be an SGFEM.

In this paper, we consider an ``interface problem" modeled by a scalar second order elliptic PDE in 2-D with piecewise smooth coefficients. We will numerically investigate the accuracy, conditioning, and the robustness of the GFEM associated with various forms of enrichments used in the literature, when applied to this problem. We will especially investigate the performance of an iterative procedure to solve the underlying linear system of the GFEM, where the stopping criterion is based on computed \textit{discretization error} and \textit{truncation error}.

In particular, we will consider GFEM with (i) the ``topological enrichment'' where a minimal order of enriched nodes are used, (ii) M-GFEM which is a generalization of topological enrichment, (iii) the ``geometrical enrichment" where all the nodes in a fixed neighborhood of the interface are enriched, and (iv) the so-called SGFEM obtained by the simple modification of M-GFEM enrichments, as suggested in \cite{BabBan}. Through numerical experiments, we show that the SGFEM is accurate for interface problems and that it does not lose accuracy as suggested in \cite{GupDuaEtAl,GupDua2EtAl,SauFri2}; thus no further modification of the enrichment is required for the interface problem (in contrast with \cite{GupDuaEtAl,GupDua2EtAl} where it was required to restore accuracy). We also show that among all the enrichments considered in the paper, the SGFEM is the only method that is well-conditioned and robust with respect to the mesh for the interface problem. These properties of the SGFEM will be proved mathematically in a forthcoming paper. One of the most important features of the current paper is the study of the performance of an iterative procedure to solve the linear system associated with M-GFEM and SGFEM. For a given error tolerance $\tau$, we have computed the solutions of the linear systems for the M-GFEM and SGFEM using the iterative method with the stopping criterion based jointly on discretization--truncation errors, as mentioned before. We observed that SGFEM requires fewer iterations and less ``wall-clock'' time than the M-GFEM.

The outline of the paper is as follows: in Section 2 we define the interface problem. In Section 3 we describe the GFEMs with various enrichments together with their convergence and conditioning properties in 1-D to communicate the idea in a simpler setting. In Section 4 we consider a ``straight interface'' problem in 2-D with no singularity, which could be viewed as a ``laboratory problem.'' We describe the GFEM with various enrichments for the straight interface problem, define an ``angle-condition'' that dictates the conditioning of the GFEM, and present various numerical results addressing the accuracy, conditioning, as well as the relation between the angle condition and conditioning. In this section we also discuss the notion of robustness with respect to the mesh. The numerical results clearly indicate that the GFEM with modified M-GFEM enrichments is indeed an SGFEM and is the most robust of all the methods considered. In Section 5 we present similar numerical results for a circular interface problem and come to the same conclusions as in Section 4. In Section 6 we describe the iterative method and the stopping criteria. We study its performance on linear systems associated with FEM, M-GFEM, and SGFEM.

\section{Formulation of the interface problem}\label{sec2}

Let $\Omega \subset \mathbb R^2$ be a bounded, simply connected domain with smooth boundary $\partial \Omega$. Consider another simply connected domain $\Omega_1 \subset \Omega$ and set $\Omega_0:= \Omega \backslash \ov{\Omega}_1$. $\Gamma:= \ov{\Omega}_0 \cap \ov{\Omega}_1$ is called the \textit{interface} and is assumed to be smooth. 

We are interested in the weak solution $u$ of the problem
\begin{equation}
-\nabla \cdot ( a\, \nabla u) = f + q\delta(\Gamma), \quad \mbox{in $\Omega$}, \label{PDE}
\end{equation}
with the boundary condition
\begin{eqnarray}
&&u = g_D, \quad \mbox{on } \partial \Omega \label{BCdir},\\
\mbox{or,} &&a\, \nabla u \cdot \bm{\vec{n}} = g_N, \quad \mbox{on } \partial \Omega, \label{BCneu}
\end{eqnarray}
where $0 < \beta_0 \le a(\bm{x}) \le \beta_1$, $a_i(\bm{x}):= a\big|_{\Omega_i}$, $f_i(\bm{x}):= f\big|_{\Omega_i}$ are smooth functions on $\ov{\Omega}_i$ for $i=0,1$; in particular, $a_i \in C^2(\ov{\Omega}_i)$ and $f_i \in C(\ov{\Omega}_i)$. Moreover, $\delta(\Gamma)$ is the Dirac function on $\Gamma$, $q(s)$ is smooth on $\Gamma$, $g_D(s)$ is smooth on $\partial \Omega$ and $g_N(s)$ is smooth on $\partial \Omega \cap \ov{\Omega}_i$; in particular, $g_D(s) \in C^2(\partial \Omega)$, $g_N(s) \in C^1(\partial \Omega \cap \ov{\Omega}_i)$, and $q(s) \in C^1(\Gamma)$, where $s$ is the arc length parameter.

The weak solution $u \in H^1(\Omega),\ u\big|_{\partial \Omega} = g_D$ of the Dirichlet boundary value problem \eqref{PDE},\eqref{BCdir} satisfies
\[
B(u,v):= \int_\Omega a\, \nabla u \cdot \nabla v\, d\bm{x} = \int_\Omega fv\, d\bm{x} + \int_\Gamma qv\, ds,\quad \mbox{for all } v \in H_0^1(\Omega).
\]
The solution of the above variational problem exists and is unique.

In the case of Neumann boundary value problem \eqref{PDE},\eqref{BCneu}, the weak solution $u \in H^1(\Omega)$ satisfies
\begin{equation}
B(u,v) = \int_\Omega fv\, d\bm{x} + \int_{\partial \Omega} g_N v\, ds + \int_\Gamma qv\, ds,\quad \mbox{for all } v \in H^1(\Omega). \label{WFneu}
\end{equation}
The solution of the above variational problem exists and is unique up to an additive constant, provided the data $f,\, g_N,\, q$ satisfy the compatibility condition
\begin{equation}
\int_\Omega f\, d\bm{x} + \int_{\partial \Omega} g_N \, ds + \int_\Gamma q\, ds = 0. \label{CompCond}
\end{equation}

Under the assumed smoothness on the input data $f, q, g_D,g_N$, the solution $u$ of the variational problem for \eqref{PDE}--\eqref{BCneu} is continuous on $\ov{\Omega}$ and $u \in C^2(\ov{\Omega}_i)$ for $i=0,1$. In particular, the solution does not have any singularity anywhere in $\ov{\Omega}$. We mention that in this paper we do not address the minimum regularity requirements on the input data for the analysis of GFEM. 

\begin{remark} Note that the data $q(s)$ appears in the right hand side of the variational problem and has no effect on the choice of enrichments in the GFEM and on the features of GFEM that we study in this paper. Thus, we will use $q(s)=0$ in all our numerical experiments.
\end{remark}

In the numerical experiments presented in this paper, we will consider $\ov{\Omega}_1$ to be a closed disk, the interface $\Gamma$ is thus a circle, away from the boundary $\partial \Omega$ (see Figure \ref{fig:circlepblm}). We will also consider a different interface problem, namely, the ``straight interface problem," where the interface $\Gamma$ is a straight line intersecting the boundary $\partial \Omega$ at two points (see Figure \ref{fig:straightpblm}). In general, the solution of this problem will have singularities at the points $\Gamma \cap \partial \Omega$. However, we will consider a manufactured solution of the straight interface problem that does not have any singularities, mimicking the property of the solution of a closed interface with data described before. We chose the straight interface problem to highlight in an easy and efficient manner the robustness properties of various GFEMs that we consider in this paper. Note that the results related to the straight interface problem are directly relevant to the problems with non-circular interfaces, where part of $\partial \Omega_1$ is a straight line. Moreover, we will consider only the Neumann problem \eqref{WFneu} in all our numerical experiments.

\section{Various GFEMs on a 1-D problem}\label{sec3}

The goal of this paper is to study certain features of various GFEMs, e.g., the accuracy and conditioning when applied to an interface problem. We describe these methods and the associated features for a 1-D problem, which will allow us to communicate the main ideas in a simpler setting.

Let $\Omega = (0,1)$, $\Omega_0 = (0,\gamma)$, $\Omega_1 = (\gamma, 1)$, where $0 < \gamma < 1$ is the interface. We consider the boundary value problem
\begin{eqnarray*}
&&- (au^\prime)^\prime = f,\quad \mbox{in }\Omega, \\
&&u(0) = 0,\ au^\prime(1) = g,
\end{eqnarray*}
where $f=1$, $g=2$ and 
$a(x)\big|_{\Omega_0} = a_0$, $a(x)\big|_{\Omega_1} = a_1$ with $a_0,a_1$ as strictly positive constants. 

The weak formulation of the above problem is
\begin{eqnarray}
&&u \in \mc{E}:= \{u \in H^1(\Omega)\, :\, u(0) = 0 \}, \nonumber \\
&& B(u,v) = F(v),\quad \mbox{for all } v \in \mc{E}, \label{InterfaceProb1d}
\end{eqnarray}
where
\[
B(u,v):= \int_\Omega a\, u^\prime\, v^\prime \, dx \mbox{ \ and \ } F(v):= \int_\Omega v\, dx + 2v(1).
\]
We denote the energy norm of $v \in \mc{E}$ by $\|v\|_{\mc{E}}:= B(v,v)^{1/2}$. 

The GFEM to approximate the solution of the above problem is a generalization of the standard FEM. We will describe various GFEMs below and state the results associated with their accuracy and conditioning. 
This section, associated with a 1-D problem, can be viewed as a conceptual synopsis of the results presented in this paper in higher dimensional interface problem.

\medskip

\noindent \textbf{GFEM:} Let $\mc{T}_h$ be the uniform mesh on $\Omega$ with nodes $x_i^h = ih$, $i \in \mc{N}^h:= \{0,1,\dots,m\}$ and elements $\tau_i^h = [x_{i-1}^h,x_i^h],\ i \in \mc{N}^h_d:= \{1,2,\dots, m\}$, where $h = 1/m$. With each node $x_i^h,\ i=1,2,\dots,m-1$, we associate the \textit{patch} $\omega_i^h = (x_{i-1}^h,x_{i+1}^h)$; for $i=0,m$, we set $\omega_0^h = (x_0^h,x_1^h)$ and $\omega_m^h = (x_{m-1}^h,x_m^h)$. Clearly, $\cup_{i=0}^m \omega_i^h = \Omega$.

Let $N_i^h$ be the usual piecewise linear ``hat-functions'' associated with the node $x_i^h$ with supp$\{N_i^h\}= \ov{\omega}_i^h$ and $N_i^h(x_i^h) = 1$. Note that the interface $\gamma \in \tau_{c}^h$ for some $c$ depending on $h$. The GFEM solution $u_h \in S^h$ satisfies
\begin{equation}
B(u_h,v) = F(v), \quad \mbox{for all }v \in S^h, \label{GFEM}
\end{equation}
where the approximation space $S^h$ is given by
\begin{equation}
S^h = S_{FEM}^h \oplus S_{ENR}^h = \{v = v_1 + v_2\, :\, v_1 \in S_{FEM}^h,\ v_2 \in S_{ENR}^h\}, \label{GFEMspace} 
\end{equation}
where
\begin{eqnarray}
&& S_{FEM}^h = \mbox{span}\{N_i^h,\, i \in \mc{N}^h_d\},\nonumber \\
\mbox{and}&&S_{ENR}^h = \mbox{span}\{wN_i^h, i \in \mc{R}^h \subset \mc{N}^h \label{S2} \}. 
\end{eqnarray}
The function $w$ in \eqref{S2} is called the \textit{enrichment} and $S_{ENR}^h$ is called the \textit{enrichment space} of the GFEM. The function $w$ is chosen such that it mimics the exact solution. The set $\{x_i^h\}_{i\in \mc{R}^h}$ denotes the set of \textit{enriched nodes}. Different GFEMs are defined by different choices of $w$ and the choice of $\mc{R}^h$. $S_{FEM}^h$ is the FE space of piecewise linear functions associated with the ``triangulation" $\mc{T}_h$, satisfying the homogeneous Dirichlet boundary condition at $x_0^h$. It is clear that dim$(S_{FEM}^h)=m$. Note that for $w \equiv 0$, we do not consider $S_{ENR}^h$ in the definition of $S^h$; thus we have $S^h = S_{FEM}^h$ and the GFEM is just the standard FEM. GFEM is thus a generalization or extension of FEM; the approximation space $S^h$ of GFEM is the standard FE space $S_{FEM}^h$ augmented with the enrichment space $S_{ENR}^h$, as given in \eqref{GFEMspace}. In the rest of the section, we will use $x_i,\tau_i,\omega_i, N_i$ for $x_i^h,\tau_i^h,\omega_i^h, N_i^h$ with an understanding that they depend on $h$.

The GFEM solution $u_h \in S^h$ is obtained in the form $u_h = \sum_{i\in \mc{N}^h_d} c_{1,i}N_i + \sum_{k\in \mc{R}^h}c_{2,k}wN_k$ by solving the linear system
\begin{equation}
\hat{\bm{A}} \hat{\bm{c}} = \hat{\bm{f}}, \label{badLinSys}
\end{equation}
where
\[
\hat{\bm{A}} = \left[ \begin{array}{cc} \hat{\bm{A}}_{11} & \hat{\bm{A}}_{12} \\ \hat{\bm{A}}_{12}^T & \hat{\bm{A}}_{22} \end{array} \right],\ \hat{\bm{c}} = \left[ \begin{array}{c} \hat{\bm{c}}_1 \\ \hat{\bm{c}}_2 \end{array} \right],\ \hat{\bm{f}} = \left[ \begin{array}{c} \hat{\bm{f}}_1 \\ \hat{\bm{f}}_2 \end{array} \right],
\]
and
\begin{eqnarray*}
&&\hat{\bm{A}}_{11} = [B(N_j,N_i)]_{i,j\in \mc{N}^h_d}, \ \ \hat{\bm{A}}_{22} = [B(wN_k,wN_l)]_{k,l \in \mc{R}^h},\\
&& \hat{\bm{A}}_{12} = [B(wN_k,N_i)]_{i \in \mc{N}^h_d\,;\, k \in \mc{R}^h},\ \ \hat{\bm{f}}_1 = [F(N_i)]_{i\in \mc{N}^h_d},\ \ \hat{\bm{f}}_2 = [F(wN_l)]_{l \in \mc{R}^h}, \\
&& \hat{\bm{c}}_1 =[c_{1,i}]_{i\in \mc{N}^h_d},\ \ \hat{\bm{c}}_2 =[c_{2,k}]_{k\in \mc{R}^h}.
\end{eqnarray*}
Note that the matrix $\hat{\bm{A}}_{11}$ is the standard FEM stiffness matrix. The matrices $\hat{\bm{A}}_{12}$ and $\hat{\bm{A}}_{22}$ depend on the enrichment $w$. Also dim$(\hat{\bm{A}}_{22}) = \mbox{card}(\mc{R}^h)$. We further note that depending on $w$ the diagonal elements of $\hat{\bm{A}}_{22}$ could be small and therefore, instead of solving the linear system \eqref{badLinSys}, we solve the linear system
\begin{equation}
\bm{A}\bm{x} = \bm{f}, \label{goodLinSys}
\end{equation}
with
\[
\bm{A} = \bm{D} \hat{\bm{A}} \bm{D} = \left[ \begin{array}{cc} \bm{A}_{11} & \bm{A}_{12}\\ \bm{A}_{12}^T & \bm{A}_{22} \end{array} \right], \ \bm{f} = \bm{D} \hat{\bm{f}},\ \bm{x} = \bm{D}^{-1}\hat{\bm{c}},
\]
where $\bm{D}$ is a diagonal matrix such that $\bm{A}$ has unit diagonal elements.

We will now describe examples of GFEM for the interface problem, where the interface $\gamma \in \mathring{\tau}_c=(x_{c-1},x_c)$; note that $\mathring{\tau}_c$ is the interior of the (closed) element $\tau_c = [x_{c-1},x_c]$.

\medskip

\textbf{Geometric GFEM:} The enrichment part $S_{ENR}^h$ is defined with the enrichment $w =\left|x-\gamma\right|$ and $\mc{R}^h = \{i\in \mc{N}^h\, :\, \left|x_i - \gamma\right| \le R\}$, where $R>0$ is fixed and independent of $h$. Then the card$\{x_i\}_{i \in \mc{R}^h} = O(h^{-1})$, and consequently, the dimension of $\hat{\bm{A}}_{22}$ is $O(h^{-1})$. The idea of using $\mathcal{R}^h$ with $R$ fixed and independent of $h$ was used in \cite{BecMinEtAl, NicRenard, SukChopEtAl}; for other references see \cite{BelGraVen,FriBel}.

\medskip

\textbf{Topological GFEM:} The enrichment part $S_{ENR}^h$ is defined with the enrichment $w = \left|x-\gamma\right|$ as before, however, $\mc{R}^h = \{x_{c-1},x_c\}$. Note that $\mc{R}^h$ is the union of nodes of the element $\tau_c$ containing the interface $\gamma$. Consequently, $\hat{\bm{A}}_{22}$ is a $2 \times 2$ matrix and the associated stiffness matrix is smaller than that of the Geometric GFEM. This idea was first used in \cite{BelyBla} in the context of crack propagation where the nodes close to the crack-tip were associated with $\mathcal{R}^h$. For other references see \cite{BelGraVen,FriBel}.

\medskip

\begin{figure}[!ht]
\centerline{\includegraphics[scale = 0.5]{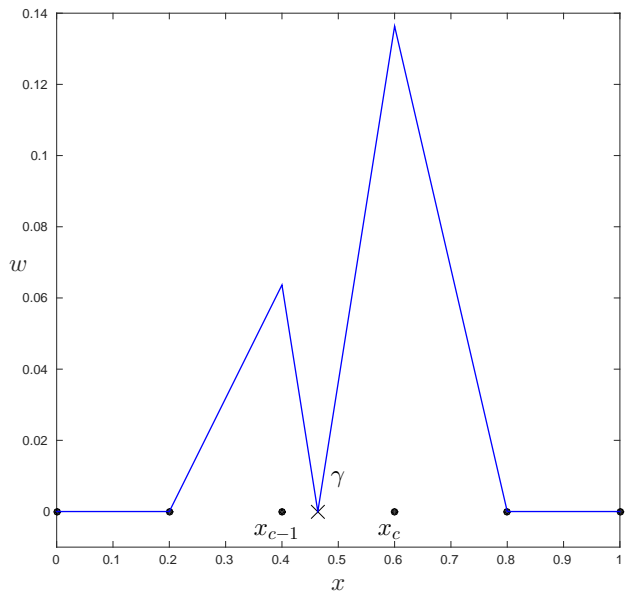}}
\caption{Enrichment function used for M-GFEM.}
\label{fig:mfct1d}
\end{figure}

\textbf{M-GFEM:} Let $w^*= \left|x-\gamma\right|$. We consider $S_{ENR}^h$ with the enrichment function $w$ that is continuous in $\Omega$, $w = w^*$ in the element $\tau_c = [x_{c-1},x_c]$ containing $\gamma$. It is \textit{linear in every element other than $\tau_c$}, and $w(x_i)=0$ for all the nodes $x_i$ \textit{except} $x_{c-1}$ and $x_c$. The graph of $w$ is of the form ``M'', as shown in Figure \ref{fig:mfct1d}, where the parameters are $\gamma = \left[2 + 1/\pi\right]h$ and $h=1/5$. We consider $\mc{R}^h = \{x_{c-2},x_{c-1},x_c,x_{c+1}\}$, which is the union of all the nodes in $\ov{\omega}_{c-1} \cup \ov{\omega}_c$; $\ov{\omega}_{c-1}, \ov{\omega}_c$ are the closures of the patches containing the interface $\gamma$. Thus $\hat{\bm{A}}_{22}$ is a $4 \times 4$ matrix and the associated stiffness matrix is smaller than that of the Geometric GFEM, but slightly bigger than that of the Topological GFEM.

The Topological GFEM and M-GFEM are local in the sense that only nodes close to the interface $x=\gamma$ are enriched and that their cardinality does not change with $h$. We mention that another local GFEM was introduced in \cite{Fries2}, referred to as the Corrected XFEM that employed a ``Ramp cutoff function." However the use of ramp-function yields the enrichment function $w$, which is quadratic in the elements $[x_{c-2},x_{c-1}], \, [x_{c},x_{c+1}]$; the enrichment in M-GFEM is piecewise linear in every element. In contrast, the Geometric GFEM is global in the sense that card$\{x_i\}_{i \in \mc{R}^h} = O(h^{-1})$; for $R$ large enough, we actually have $\mc{R}^h = \mc{N}^h$, i.e., every node of FE mesh $\mc{T}_h$ could be enriched in the Geometric GFEM for a large choice of $R$. 

We further note that if the interface $x=\gamma$ is at one of the nodes $x_i$, then no enrichment is used, i.e., standard FEM can be used. However, if $x=\gamma$ is close to a node, the round-off error may create a serious problem in GFEM with the enrichments described above. In such a situation there has to be a safety check and it is advisable not to use the enrichment; we only use $S^h_{FEM}$. We note that not using the $S^h_{ENR}$ part may result in loss of accuracy in the approximation, but this loss is of much smaller scale compared to the round-off error that may result if we used $S^h_{ENR}$. For M-GFEM, we have observed that not using any enrichment when $\min \left\lbrace \gamma - x_{c-1},x_c - \gamma \right\rbrace\le 10^{-14}h$ does not affect the accuracy of the computed solution $u_h$. The factor $10^{-14}$ in the safety-check depends on the machine precision of the computer.

First we highlight the performance of various GFEMs with respect to their accuracy. Let $\epsilon^h:= \|u-u_h\|_\mc{E}$, where $u_h$ is the computed solution using one of the GFEMs. For the Topological GFEM, one can show that $\epsilon^h \le Ch^{1/2}$ -- similar to the well-known result for FEM with uniform mesh and $\gamma \in \mathring{\tau}_c$. But for the Geometric GFEM and the M-GFEM, the rate of convergence is higher, namely, $\epsilon^h \le Ch$. The order of convergence for Geometric GFEM and M-GFEM applied to an interface problem is thus the same as the one for the FEM applied to a problem with smooth solution.

However, the conditioning of a GFEM could be much worse than the conditioning of the FEM. Consequently, solving the linear system \eqref{goodLinSys} could be extremely difficult. In fact, the condition number $\kappa_2(\bm{A})$ for the GFEM depends on the ``angle'' between the spaces $S_{FEM}^h$ and $S_{ENR}^h$, which will be defined precisely in the next section. It has been shown in \cite{BabBan,QzhangBabBan} that if the angles between the spaces $S_{FEM}^h$ and $S_{ENR}^h$ are ``not too small,'' then $\kappa_2(\bm{A})= O(h^{-2})=\kappa_2(\bm{A}_{11})$, i.e., the conditioning of the GFEM is not worse than that of the standard FEM. It was shown in \cite{BabBan} that the ``angle'' could become very small for typical enrichments of GFEM used in practice, and $\kappa_2(\bm{A})$ could be $O(h^{-4})$.

Therefore, to design a well-conditioned GFEM, one has to choose the enrichment function $w$ and the enrichment space $S_{ENR}^h$ such that the ``angle" between $S_{FEM}^h$ and $S_{ENR}^h$ ``is not too small'', i.e., stays bounded away from 0.

\medskip

\textbf{Stable GFEM (SGFEM):} A GFEM is called an SGFEM if (a) $\epsilon^h \le Ch$ and (b) ``angle'' between $S_{FEM}^h$ and $S_{ENR}^h$ stays bounded away from $0$ for all $h$. Specifically, for the interface problem \eqref{InterfaceProb1d}, we let $w^*= \left|x-\gamma\right|$. $S_{ENR}^h$ is defined with the enrichment $w = w^* - \ih w^*$, where $\ih w^*$ is the piecewise linear interpolant of $w^*$ with respect to the triangulation $\mc{T}_h$. Clearly $w$ is continuous in $\Omega$ and $w=0$ outside $\tau_c$. Again we consider $\mc{R}^h = \{x_{c-1},x_c\}$, and consequently $\hat{\bm{A}}_{22}$ is a $2 \times 2$ matrix. The GFEM with enrichment $w$ defined above yields $\epsilon^h \le Ch$, the angle between the associated $S_{ENR}^h$ and $S_{FEM}^h$ is bounded away from $0$, and $\kappa_2(\bm{A}) =O(h^{-2})$. Thus this GFEM is indeed an SGFEM. It is local, and \textit{will be referred to as the SGFEM in this paper}. This enrichment and the associated GFEM was used in \cite{BabBan,MoeCloEtAl}. Note that the same procedure applied to the enrichment in Corrected XFEM will yield $\mathcal{R}^h = \{ x_{c-2},x_{c-1},x_{c},x_{c+1} \}$, i.e., Corrected XFEM will require more degrees of freedom.

\medskip

We mention that M-GFEM is well-conditioned. But in higher dimensions, the conditioning of M-GFEM is not robust with respect to the position of the interface $\Gamma$ with respect to the mesh, which we will show in the next section. On the other hand, the Geometric GFEM is badly conditioned in the sense that $\kappa_2(\bm{A}) = O(h^{-4})$. Note however that one has to choose $h$ small enough, depending on $R$, to see this effect.

In summary, we say that the Topological GFEM is not as accurate as other GFEMs. The Geometric GFEM, though accurate, is not well conditioned for all $h$. The M-GFEM is accurate but the conditioning is not robust in higher dimensions. However, the SGFEM is accurate as well as robustly well-conditioned. These features will be shown for 2-D problems in the later sections of this paper.

\section{Straight interface problem}\label{sec4}

In this section, we discuss the GFEM applied to a 2-D problem with a straight interface. We consider a specific manufactured problem such that the solution does not have any singularities. This problem will allow us to easily show the process of extending the 1-D ideas, presented in Section 3, to 2-D problems without the technicalities involved in a general interface problem. This problem will allow us to compare the robustness of various GFEMs considered in this paper. Moreover the straight interface problem could actually be viewed as a ``laboratory problem.''

\begin{figure}[!ht]
\centerline{\includegraphics[scale = 0.65]{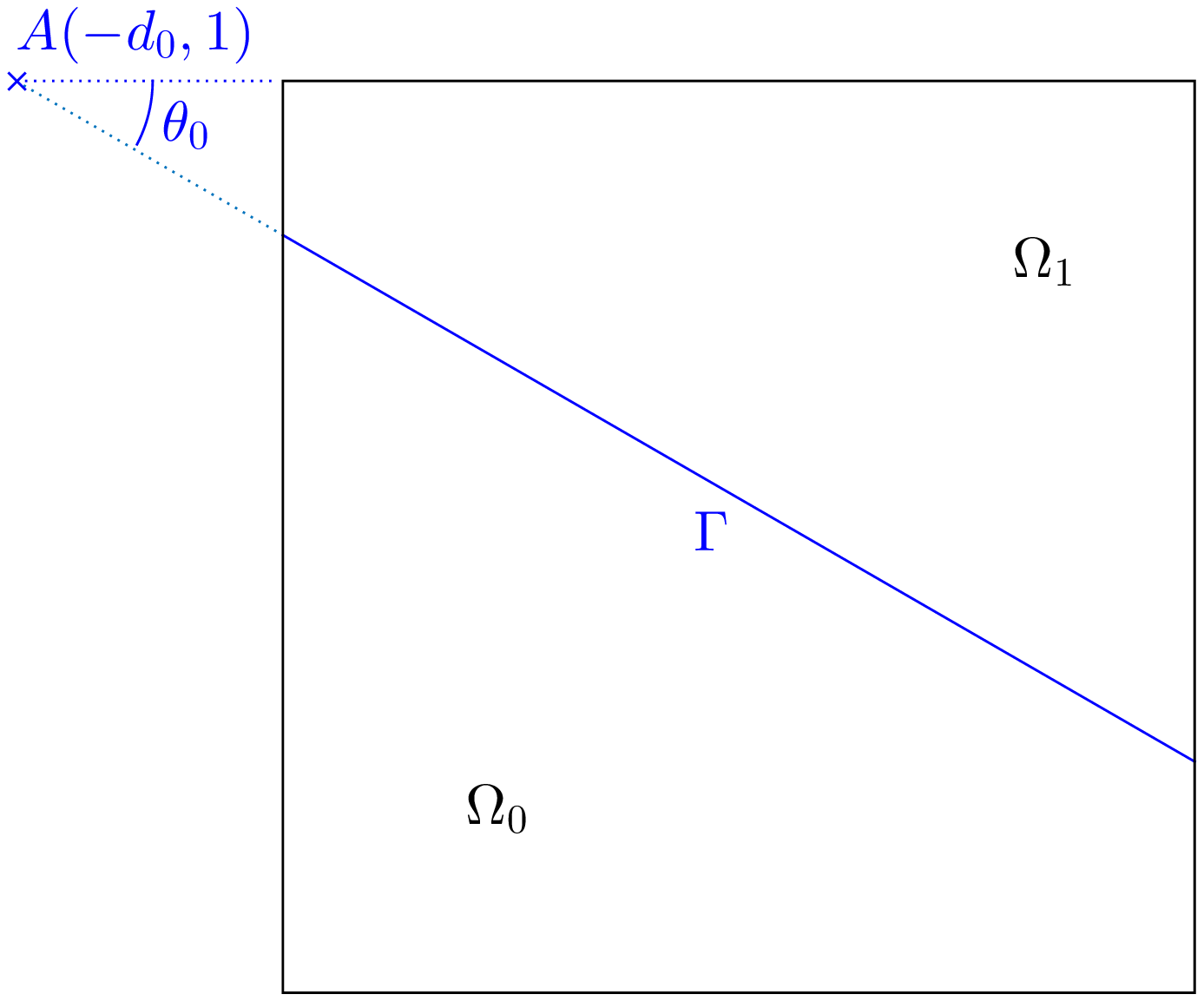}}
\caption{Straight interface problem.}
\label{fig:straightpblm}
\end{figure}

Consider the domain $\Omega = (0,1) \times (0,1)$. For a given $d_0>0$ and $0 < \theta_0 < \tan^{-1}(1/d_0)$, let $\Gamma:= \{\bm{x} \in \Omega\,:\, \gamma(\bm{x})=0\}$ be the straight interface, where $\gamma(\bm{x})=0$ is the straight-line passing through the point $A(-d_0,1)$ with slope $\mc{M}=-\tan(\theta_0)$, as shown in Figure \ref{fig:straightpblm}, where we have chosen $d_0 = 1-1/\sqrt{2}$ and $\theta_0 = \pi/6$. We set $\Omega_0:= \Omega \cap \{\bm{x}\,:\, \gamma(\bm{x}) <0\}$ and $\Omega_1:= \Omega \cap \{\bm{x}\,:\, \gamma(\bm{x}) >0\}$. Note that varying $d_0$ and $\theta_0$, we can vary the interface $\Gamma$.

We consider the problem \eqref{WFneu} with $f\equiv 0$, $q \equiv 0$, and $g_N$ satisfying the compatibility condition \eqref{CompCond}. Moreover, we consider $a_i(\bm{x}) = a_i$, $i=0,1$, where $a_0, a_1$ are strictly positive constants. The solution $u \in \mc{E} = H^1(\Omega)$ of \eqref{WFneu} exists and is unique up to an additive constant. We set $\|u\|_{\mc{E}}:= B(u,u)^{1/2}$.

For the computations presented in this section, we will consider the manufactured solution of \eqref{WFneu}, given by
\begin{equation}
u_{ex} = \left\{ \begin{array}{lr}
A_0r^\alpha \cos[\alpha(\theta - \theta_0)] + B_0 r^\alpha \sin[\alpha(\theta - \theta_0)] + C,& \theta \le \theta_0, \\
A_1r^\alpha \cos[\alpha(\theta - \theta_0)] + B_1 r^\alpha \sin[\alpha(\theta - \theta_0)] + C,& \theta \ge \theta_0, \end{array} \right. \label{ManufacSol}
\end{equation}
where $(r,\theta)$ is the polar coordinate centered at $A$ and the polar line $\{\bm{x}=(x_1,x_2): -d_0 < x_1, x_2=1\}$. We choose $A_0,A_1,B_0,B_1$ such that $u_{ex}$ is continuous in $\Omega$ and $a(\bm{x})\, \frac{\partial u_{ex}}{\partial n}$ is continuous across the interface $\Gamma\, (\theta = \theta_0)$. Then, we choose $C$ such that $u_{ex}(0,0) = 0$. We also consider $\alpha \ge 1$. Clearly, $u_{ex}$ is continuous in $\Omega$ with no singularity in $\ov{\Omega}$. Also $f = 0$ and $g_N$, in \eqref{WFneu}, is obtained from $u_{ex}$ using \eqref{BCneu}.


We now describe the GFEM in 2-D. Let $\mc{T}_h$ be a uniform finite element triangulation of $\Omega$ with nodes $\bm{x}_{\bi} = (i_1h,i_2h)$, where for a given positive integer $m$, we define $h = \frac{1}{m}$ and $\bi \in \mc{N}^h:= \{(i_1,i_2)\,:\, i_1,i_2 = 0,1,2,\cdots,m\}$. We denote the set of elements $\tau$, which are closed triangles with nodes as their vertices, by $E$. For each node $\bm{x_i}$, we define $\ov{\omega}_{\bi} = \{\cup \tau\,: \, \bm{x_i} \mbox{ is a vertex of } \tau \in E\}$. For the given triangulation $\mc{T}_h$, $\ov{\omega}_{\bi}$ is the union of $1,2,3$ or 6 elements with the vertex $\bm{x_i}$ depending on its position in $\Omega$. The open set $\omega_{\bi}$ is the patch associated with the node $\bm{x_i}$. It is clear that $\Omega = \cup_{\bi \in \mc{N}^h} \omega_{\bi}$. Let $N_{\bi}$ be the usual piecewise linear hat-function associated with node $\bm{x_i}$ with supp$\{N_{\bi}\} = \ov{\omega}_{\bi}$ and $N_{\bi}(\bm{x_i}) = 1$. We set
\[
E_\Gamma:= \{\tau \in E \,:\, \mathring{\tau} \cap \Gamma \ne \emptyset \}.
\]

The approximation space $S^h$ of the GFEM in 2-D is constructed similarly to the description given in Section \ref{sec3}. We set
\[
S_{FEM}^h = \mbox{span}\{N_{\bi},\, \bi \in \mc{N}_d^h\}, \quad \mbox{where } \mc{N}_d^h = \{\bi \in \mc{N}^h\,:\, \bi \ne (0,0)\}.
\]
The functions in the space $S_{FEM}^h$ vanish at the node $(0,0)$ and dim$\{S_{FEM}^h\} = (m+1)^2 -1$. For a given enrichment function $w$ that mimics the exact solution, we also define the \textit{enrichment space}
\[
S_{ENR}^h = \mbox{span}\{wN_{\bi},\ \bi \in \mc{R}^h \subset \mc{N}^h \}.
\]
The particular choice of the enrichment function $w$ and of the set of indices $\mc{R}^h \subset \mc{N}^h$ defines a distinct GFEM. The set $\{\bm{x_i}\}_{\bi \in \mc{R}^h}$ denotes the set of \textit{enriched nodes}. The approximation space $S^h \subset \mc{E}$ of GFEM is given by \eqref{GFEMspace}, namely, $S^h = S_{FEM}^h \oplus S_{ENR}^h$.

The GFEM solution $u_h \in S^h$ satisfies the finite dimensional problem \eqref{GFEM} with $F(v) = \int_{\partial \Omega} g_Nv\, dx$ and is given in the form $u_h = \sum_{\bi \in \mc{N}_d^h} c_{1,\bi} N_{\bi} + \sum_{\bm{k} \in \mc{R}^h} c_{2,\bm{k}}wN_{\bm{k}}$, where $c_{1,\bi},\, c_{2,\bm{k}}$ is the solution of the linear system \eqref{badLinSys}. Note that we solve the diagonally scaled linear system \eqref{goodLinSys}, instead of \eqref{badLinSys}. Note also that if $w\equiv 0$, i.e., if no enrichment is used, $S^h = S_{FEM}^h$ and the GFEM is the standard FEM.

We will now describe various GFEMs based on the particular choices of the enrichment function $w$ and the set of indices $\mc{R}^h$ for the enriched nodes. The enrichment function $w$ for the interface problem is based on the so-called \textit{distance function}
\begin{equation}
w^*(\bm{x}):= \mbox{dist}(\bm{x},\Gamma). \nonumber
\end{equation}
Note that $w^*(\bm{x})$ is continuous in $\Omega$, it is linear in $\Omega_0$ and $\Omega_1$, and $w^*(\bm{x})=0$ for $\bm{x} \in \Gamma$.

\medskip

\textbf{Geometric GFEM:} The enrichment space $S_{ENR}^h$ is constructed with the enrichment function $w(\bm{x}) = w^*(\bm{x})$, whereas the set of indices $\mc{R}^h$ is given by $\mc{R}^h=\{\bi \in \mc{N}^h\,:\, \mbox{dist}(\bm{x_i}, \Gamma) \le R\}$ for a fixed $R$, independent of $h$. Note that unlike in Section \ref{sec3}, card$\{\bm{x_i}\}_{\bi \in \mc{R}^h} = O(h^{-2})$. The set of enriched nodes thus contains all the nodes within a fixed distance $R$ from the interface.


\medskip

\textbf{Topological GFEM:} The same enrichment function $w(\bm{x})=w^*(\bm{x})$ is used to define $S_{ENR}^h$ as in the Geometric GFEM. Here however, we use 
$\mc{R}^h = \{ \bi\in \mc{N}^h\,:\, \omega_{\bi} \cap \Gamma \ne \emptyset \}$. Note that the set of enriched nodes $\{\bm{x_i}\,:\, i \in \mc{R}^h\}$ is the union of all the vertices of the elements $\tau \in E_\Gamma$. Again unlike in Section \ref{sec3}, card$\{\bm{x_i}\}_{\bi \in \mc{R}^h} = O(h^{-1})$.

\medskip

\textbf{M-GFEM:} The enrichment function used in M-GFEM is slightly different than $w$ used in Geometrical or Topological GFEM and is given by
\begin{eqnarray*}
&&w(\bm{x}) = \left\{ \begin{array}{ll}
w^*(\bm{x}), & \bm{x} \in \tau \in E_\Gamma, \\
\mbox{linear function}, & \bm{x} \in \tau \in E \backslash E_\Gamma,
\end{array} \right. \\
&& w(\bm{x_i}) = 0, \quad \bm{x_i} \mbox{ is not a vertex of } \tau \in E_\Gamma.
\end{eqnarray*}
Furthermore, the indices of enriched nodes are given by
\[
\mc{R}^h=\{\bi \in \mc{N}^h\,:\, \bm{x_i} \in \bigcup_{\omega_{\bm{j}} \cap \Gamma \ne \emptyset} \ov{\omega}_{\bm{j}} \}.
\]
 Note that unlike in Section \ref{sec3}, card$\{\bm{x_i}\}_{\bi \in \mc{R}^h} = O(h^{-1})$.
 
 \medskip
 
\begin{figure}[!ht]
\centerline{\includegraphics[scale = 0.65]{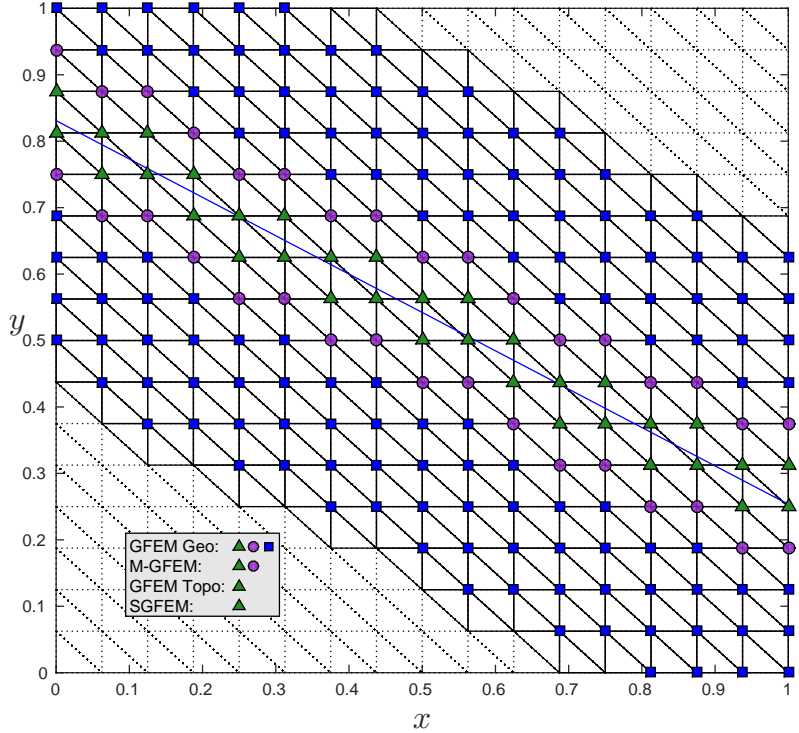}}
\caption{Nodes enriched for the straight interface problem.}
\label{fig:straightpatches}
\end{figure}

The enriched nodes $\bm{x_i},\bi \in \mc{R}^h$ are shown in Figure \ref{fig:straightpatches}, where $h=1/16$, $d_0 = 1-1/\sqrt{2}$, $\theta_0 = \pi/6$ and $R = 1/3$.

\begin{figure}[!ht]
\centerline{\includegraphics[scale = 0.65]{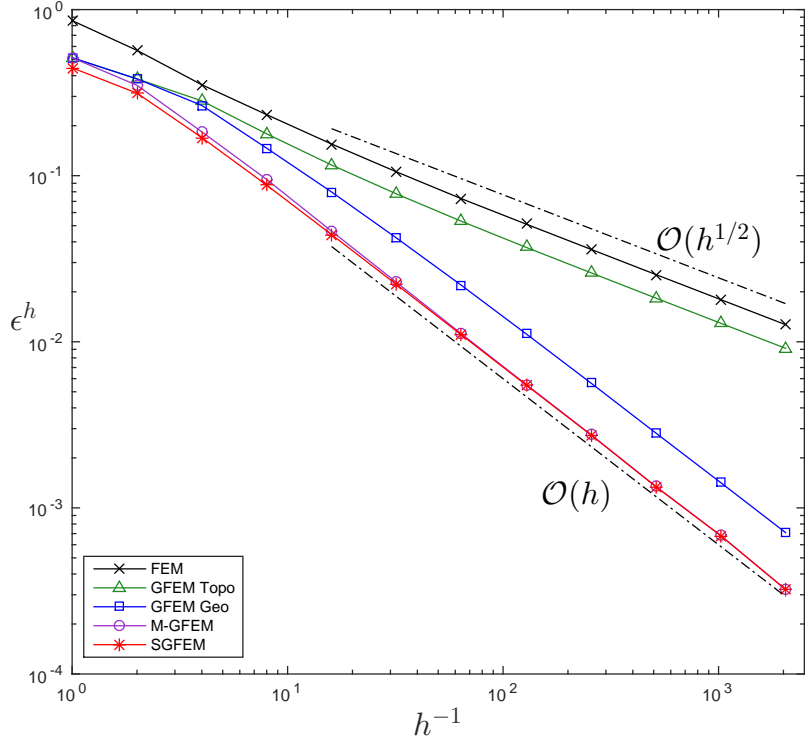}}
\caption{Relative error in the energy norm against $h$.}
\label{fig:straighterrorh}
\end{figure}

First, we consider the exact solution $u\in \mc{E}$ of \eqref{WFneu} given in \eqref{ManufacSol} with $d_0 = 1-1/\sqrt{2}$, $\theta_0 = \pi/6$, $a_0 = 1$ and $a_1 = 10$. Note that the interface $\Gamma$ is not aligned with the mesh. We computed the error $\|u-u_h\|_\mc{E}$, where $u_h$ is the solution of \eqref{GFEM} associated with the GFEMs described above (we chose $R = 1/3$ for Geometric GFEM). The log-log plot of the (relative) error is given in Figure \ref{fig:straighterrorh}. It is clear that Geometric GFEM and M-GFEM yield the convergence of $O(h)$, whereas the order of convergence for the Topological GFEM is only $O(h^{1/2})$. This suboptimal order of convergence for Topological GFEM has been reported in the literature \cite{LabPomEtAl,StaBudEtAl,SukChopEtAl}.

\begin{figure}[!ht]
\centerline{\includegraphics[scale = 0.65]{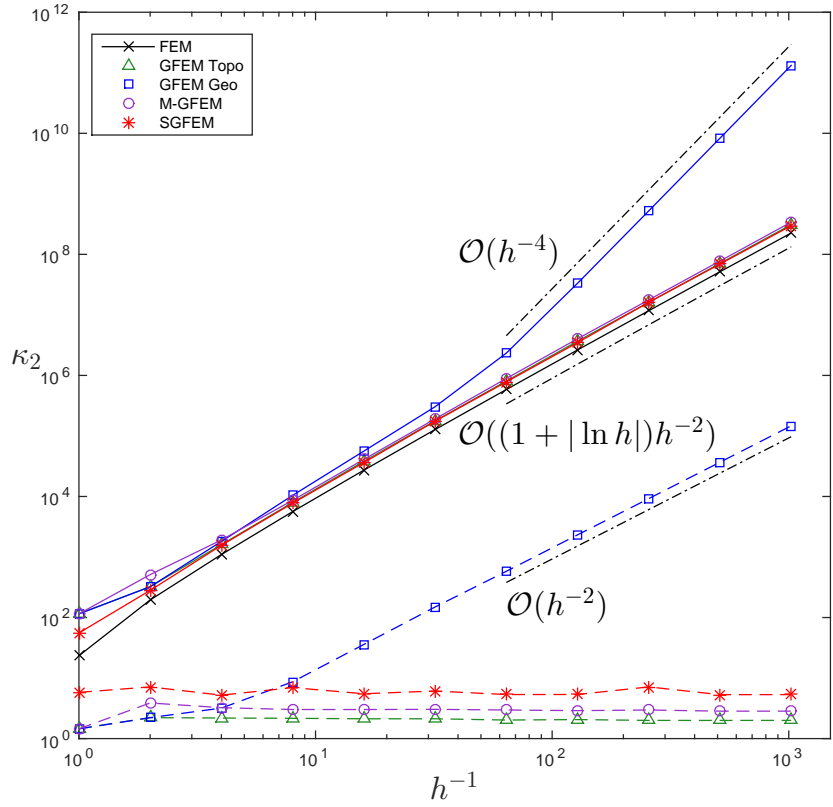}}
\caption{Condition number of the scaled stiffness matrices $\bm{A}$ (solid lines) and $\bm{A}_{22}$ (dashed lines) against $h$.}
\label{fig:straightcondh}
\end{figure}

The approximation property of GFEM, as described in the last paragraph, is undoubtedly very important. However, it is equally important that the GFEM be well-conditioned in order for the linear system \eqref{goodLinSys} to be solved efficiently. Towards this end, we first compute the condition number of the scaled stiffness matrix $\bm{A}$ in \eqref{goodLinSys} for different values of $h$. The results are collected in Figure \ref{fig:straightcondh} (solid lines), again with $d_0 = 1-1/\sqrt{2}$, $\theta_0 = \pi/6$, $R = 1/3$, $a_0 = 1$ and $a_1 = 10$. It is clear that for the Topological GFEM and M-GFEM, the condition number $\kappa_2(\bm{A}) = O(h^{-2})$. We mention that $\kappa_2(\bm{A}_{11}) = O(h^{-2})$, where $\bm{A}_{11}$ is the stiffness matrix of the standard FEM. In other words, the conditioning of the Topological GFEM and the M-GFEM is of the same order as that of a standard FEM. On the other hand, $\kappa_2(\bm{A}) = O(h^{-4})$ for the Geometrical GFEM, which is much worse than the Topological GFEM and the M-GFEM. We will show later in this paper that conditioning also plays a major role in the iterative solution of \eqref{goodLinSys}. We also mention that $\kappa_2(\bm{A}_{22})$ is bounded for all considered GFEMs, except for Geometric GFEM for which we have $\kappa_2(\bm{A}_{22})=O(h^{-2})$ (see dashed lines in Figure \ref{fig:straightcondh}).

The conditioning of GFEM depends on the ``angle'' between the spaces $S_{FEM}^h$ and $S_{ENR}^h$, which is characterized by the quantity
\[
\cos\left(\vartheta (S_{FEM}^h,S_{ENR}^h)\right):= \max_{u\in S_{FEM}^h,\, v\in S_{ENR}^h} \frac{\left|B(u,v)\right|}{\|u\|_\mc{E} \|v\|_\mc{E}},
\]
where $\vartheta (S_{FEM}^h,S_{ENR}^h) \in \left[0,90^\circ\right]$ could be interpreted as the smallest angle between $S_{FEM}^h$ and $S_{ENR}^h$ and depends on the specific enrichment function $w(\bm{x})$ used in the GFEM. In particular, if $\vartheta (S_{FEM}^h,S_{ENR}^h) = 90^\circ$, then we could write $S^h = S_{FEM}^h \overset{\perp}{\oplus} S_{ENR}^h$, i.e., the two spaces are in orthogonal direct sum. On the other hand, if $\vartheta (S_{FEM}^h,S_{ENR}^h) = 0$, then we simply have $S^h = S_{FEM}^h + S_{ENR}^h$, and the sum is no longer direct. When inbetween, $\vartheta (S_{FEM}^h,S_{ENR}^h) \in \left(0,90^\circ\right)$, we have $S^h = S_{FEM}^h \oplus S_{ENR}^h$, and the sum is direct but not orthogonal. It has been proved in \cite{BabBan,QzhangBabBan} that if the angles between the spaces $S_{FEM}^h$ and $S_{ENR}^h$ are ``not too small'', i.e., if there exist positive constants $C, C_1,C_2$ such that
\begin{equation}
\vartheta (S_{FEM}^h,S_{ENR}^h) \ge C > 0, \label{AngleCond1}
\end{equation} 
and if 
\begin{equation}
C_1 \le \kappa_2(\bm{A}_{22}) \le C_2, \label{Assump2}
\end{equation}
where $\bm{A}_{22}$ as in \eqref{goodLinSys}, then 
\begin{equation}
\kappa_2(\bm{A})= O(h^{-2})=\kappa_2(\bm{A}_{11}), \label{MainResult}
\end{equation}
i.e., the conditioning of GFEM, with $S_{FEM}^h$ and $S_{ENR}^h$ satisfying the above conditions, is not worse than that of the standard FEM. Note however that conditions \eqref{AngleCond1}--\eqref{Assump2} are sufficient conditions for \eqref{MainResult}, i.e., they guarantee the well-conditioning of the GFEM. We further note that \eqref{MainResult} holds even when the condition \eqref{Assump2} is replaced by $\kappa_2(\bm{A}_{22}) = O(h^{-2})$. Moreover, since the functions in $S_{FEM}^h$ vanish only at the node $(0,0)$, it is theoretically known \cite{Bochev2005} that $\kappa_2(\bm{A}_{11}) = O[h^{-2}(1 + \left|\ln h \right|)]$. In this paper, we will not consider the factor $\left|\ln h \right|$, and consider instead $\kappa_2(\bm{A}_{11}) = O(h^{-2})$, as done in \eqref{MainResult}. The results in \cite{BabBan,QzhangBabBan} are quite general. As long as $S^h_{ENR}$ associated with \textit{any chosen enrichment function satisfies the conditions \eqref{AngleCond1}--\eqref{Assump2}, the well-conditioning of the associated GFEM is guaranteed.} This is not only true for interface problems, but for any problem where GFEM is used. In fact, the conditions \eqref{AngleCond1}--\eqref{Assump2} may help to construct enrichments leading to well-conditioned GFEM.

\begin{figure}[!ht]
\centerline{\includegraphics[scale = 0.65]{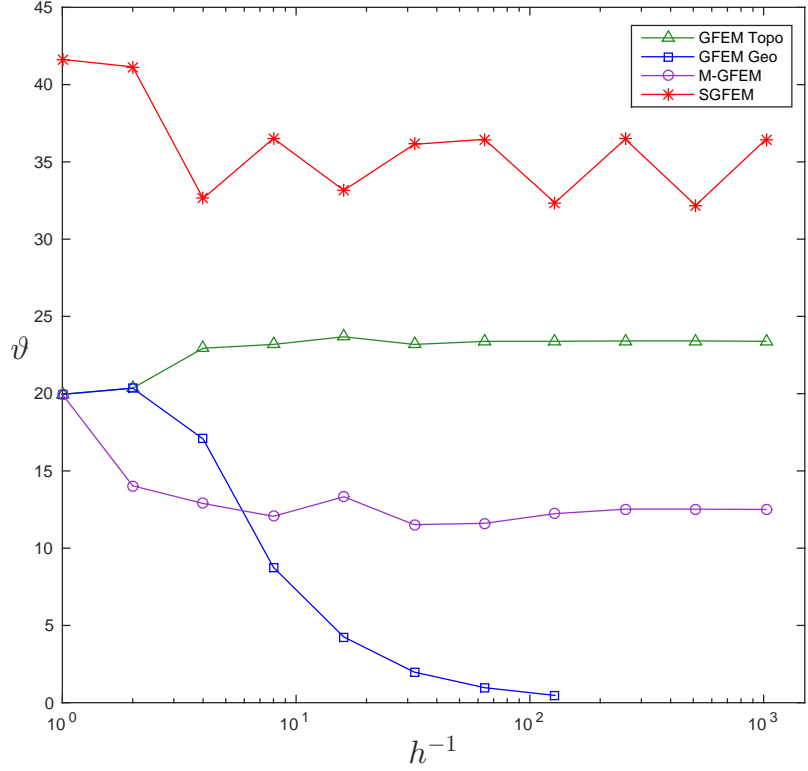}}
\caption{Angle (in degrees) between the spaces $S_{FEM}^h$ and $S_{ENR}^h$ against $h$.}
\label{fig:straightangleh}
\end{figure}

To investigate the dependence of $\kappa_2(\bm{A})$ on the angle between $S_{FEM}^h$ and $S_{ENR}^h$, we computed the angle for Geometric GFEM, Topological GFEM, and the M-GFEM for different values of $h$. The results are displayed in Figure \ref{fig:straightangleh}, again with $d_0 = 1-1/\sqrt{2}$, $\theta_0 = \pi/6$, $R = 1/3$, $a_0 = 1$ and $a_1 = 10$. The angle could be obtained by solving a generalized eigenvalue problem that we present in Appendix \ref{app:angle}. It is clear from Figure \ref{fig:straightangleh} that the angle for the Topological GFEM and M-GFEM remain bounded away from $0$ for all the values of $h$, thus illuminating the result presented above since we have seen in Figure \ref{fig:straightcondh} that $\kappa_2(\bm{A})= O(h^{-2})$ for the Topological GFEM and M-GFEM. Figure \ref{fig:straightangleh} also shows that the angle for the Geometric GFEM approaches 0 as $h$ gets smaller. Moreover, Figure \ref{fig:straightcondh} indicates that $\kappa_2(\bm{A}) = O(h^{-4}) \gg O(h^{-2})$. This suggests that the condition \eqref{AngleCond1} could be a necessary condition for \eqref{MainResult}.

Therefore, well-conditioning (of the same order as the standard FEM) of a system could be guaranteed if the GFEM uses enrichments that yield accurate approximation and if the angle between $S_{FEM}^h$ and $S_{ENR}^h$ is uniformly bounded away from $0$.

\medskip

\textbf{Stable GFEM (SGFEM):} A GFEM is called an SGFEM if (a) it yields the optimal rate of convergence, and (b) it satisfies the conditions \eqref{AngleCond1}--\eqref{Assump2}. For the straight interface problems, a special enrichment could be obtained by a simple modification of the distance function $w^*(\bm{x})$ and the associated GFEM is indeed an SGFEM. The enrichment is defined as
\begin{equation}
w(\bm{x}) = w^*(\bm{x}) - \mc{I}_hw^*(\bm{x}), \label{SGFEMenr}
\end{equation}
where $\mc{I}_hw^*(\bm{x})$ is the piecewise linear interpolant of $w^*(\bm{x})$. Since $w^*(\bm{x})$ is linear on $\Omega_0$ and $\Omega_1$, it is easy to see that 
\[
\mbox{supp}\{w(\bm{x})\} = \bigcup_{\tau \in E_\Gamma} \tau.
\]
Moreover, the indices of the enriched nodes are 
\begin{equation}
\mc{R}^h = \{\bi \in \mc{N}^h\,:\, \omega_{\bi} \cap \Gamma \ne \emptyset \}, \label{SGFEMenrNodes}
\end{equation}
which is the same as the $\mc{R}^h$ used in Topological GFEM (the enriched nodes are shown in Figure \ref{fig:straightpatches}).

In Figures \ref{fig:straighterrorh}, \ref{fig:straightcondh} and \ref{fig:straightangleh}, we have plotted the error $\|u-u_h\|_{\mc{E}}$, the condition number $\kappa_2(\bm{A})$ and the angle between $S_{FEM}^h$ and $S_{ENR}^h$, with respect to $h$, for the GFEM with enrichment given by \eqref{SGFEMenr} and $\mc{R}^h$ as in \eqref{SGFEMenrNodes}. It is clear that the method yields the optimal order of convergence, i.e., $\|u-u_h\|_{\mc{E}} = O(h)$, $\kappa_2(\bm{A}) = O(h^{-2})$, and the angle between $S_{FEM}^h$ and $S_{ENR}^h$ is bounded away from $0$ for all the values of $h$ considered in the experiment. The condition \eqref{AngleCond1} is thus satisfied. We have also checked (see Figure \ref{fig:straightcondh}, dashed lines) that \eqref{Assump2} is satisfied. Thus the GFEM with enrichment given in \eqref{SGFEMenr}--\eqref{SGFEMenrNodes} is indeed an SGFEM; \textit{we will refer to this GFEM as SGFEM in the rest of this paper}. We further note in Figure \ref{fig:straightangleh} that the angle between $S_{FEM}^h$ and $S_{ENR}^h$ for the SGFEM is larger than that of the M-GFEM -- this feature is central to the iterative solution of \eqref{goodLinSys}, which we will show later in this paper.

\begin{remark}
Note that the enrichment given in \eqref{SGFEMenr} and the set of enriched nodes indexed by $\mc{R}^h$ in \eqref{SGFEMenrNodes} was introduced in \cite{MoeCloEtAl}. However, the conditioning of the GFEM or the angle between $S_{FEM}^h$ and $S_{ENR}^h$ was not discussed there.
\end{remark}

\begin{remark}
It is important to note that modifying an enrichment by subtracting the piecewise linear interpolant, as we have done in \eqref{SGFEMenr}, may not yield an SGFEM for other problems. It has been shown in \cite{GupDuaEtAl} that for the crack propagation problems, modification of enrichment by subtracting a linear interpolant may yield inaccurate solutions. An additional modification of the enrichment function is needed to yield accurate solutions for such problems. However, the modification is certainly successful for the interface problems, as those considered in this paper.
\end{remark}

\begin{figure}[!ht]
\centerline{\includegraphics[scale = 0.65]{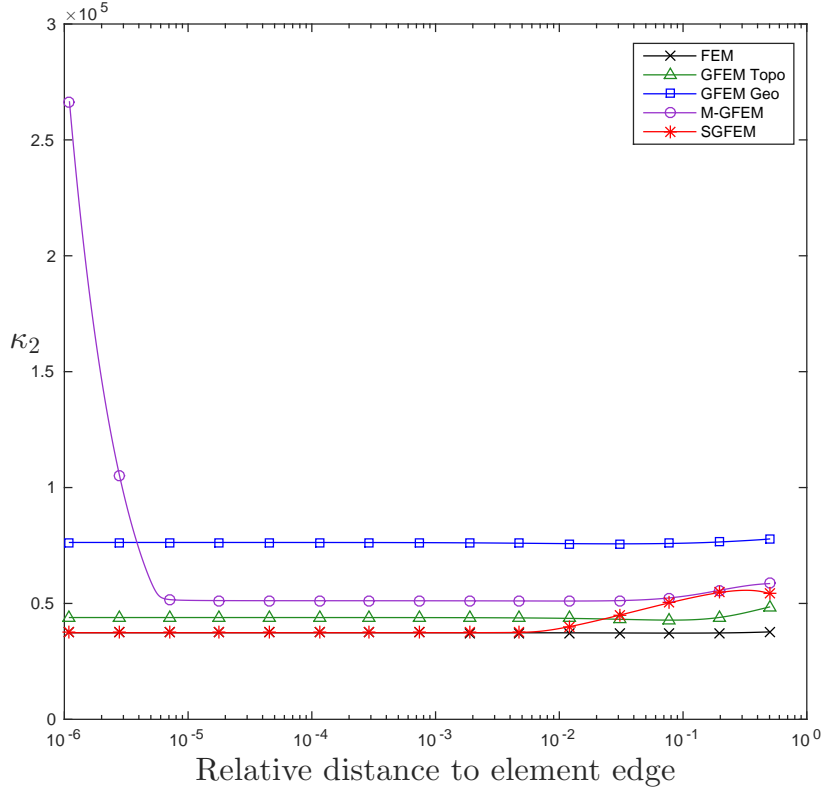}}
\caption{Condition number of the scaled stiffness matrix $\bm{A}$ when the interface gets closer to the mesh lines.}
\label{fig:straightcondd0}
\end{figure}

\begin{figure}[!ht]
\centerline{\includegraphics[scale = 0.65]{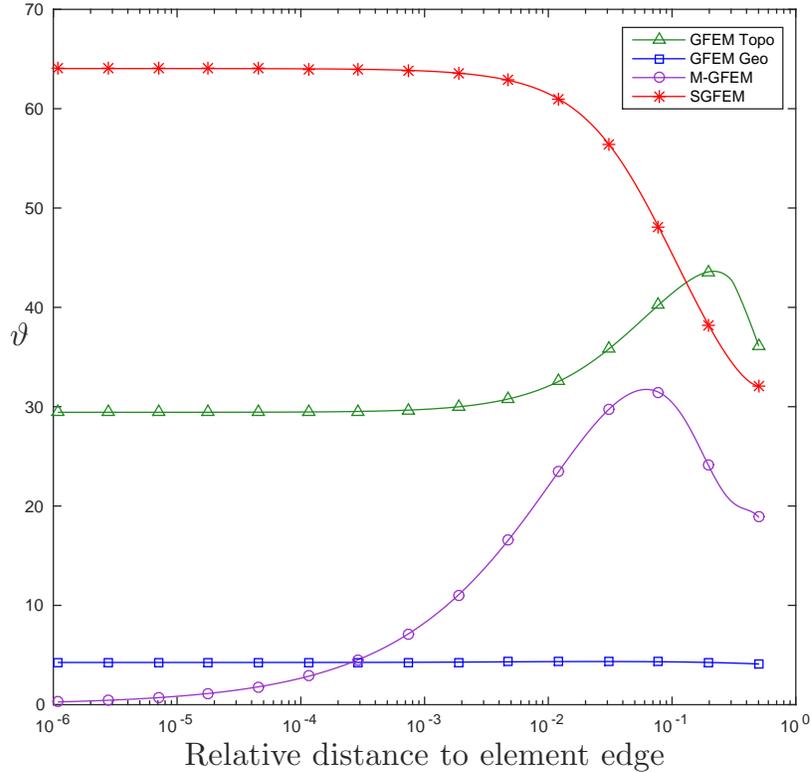}}
\caption{Angle (in degrees) between the spaces $S_{FEM}^h$ and $S_{ENR}^h$ when the interface gets closer to the mesh lines.}
\label{fig:straightangled0}
\end{figure}

It has been reported in the literature \cite{FriBel} that GFEM may become ill-conditioned if the interface $\Gamma$ is close to the mesh lines. To investigate this problem, we consider the manufactured solution \eqref{ManufacSol} with $\theta_0 = \pi/4$. In this case, the interface $\Gamma$ is parallel to some of the mesh lines associated with the triangulation $\mc{T}_h$. We control the distance of $\Gamma$ to the mesh line by controlling the parameter $d_0$ (see Figure \ref{fig:straightpblm}). We have fixed $h=1/16$ and have plotted the condition number $\kappa_2(\bm{A})$ for Topological GFEM, Geometric GFEM, M-GFEM, and SGFEM, as the interface $\Gamma$ gets closer to the mesh line in Figure \ref{fig:straightcondd0}. In this scenario, the other parameters are $R = 1/6$, $a_0 = 1$ and $a_1 = 10$. It is clear that $\kappa_2(\bm{A})$ for M-GFEM ``blows-up'' as $\Gamma$ gets closer to the mesh line; there is no appreciable change in $\kappa_2(\bm{A})$ for other GFEMs considered here. We mention that $\kappa_2(\bm{A}_{22})$ stays bounded for all GFEMs considered, even when the interface is relatively close to the mesh lines, and even for M-GFEM (not shown in Figure \ref{fig:straightcondd0}). In Figure \ref{fig:straightangled0}, we have plotted the angle between $S_{FEM}^h$ and $S_{ENR}^h$ for the fixed $h=1/16$. We clearly see that the angle for M-GFEM goes to $0$ as $\Gamma$ gets closer to the mesh line; angles for other GFEMs approach different but fixed values, bounded away from $0$. This shows that the \textit{conditioning of M-GFEM is not robust with respect to the position of the interface to the edges of the mesh}. We mention that in a forthcoming paper, we will prove that the GFEM with enrichment \eqref{SGFEMenr} satisfies the conditions \eqref{AngleCond1}--\eqref{Assump2}, where the constants are independent of $h$ and of the position of $\Gamma$. In other words, \textit{SGFEM is robust with respect to the position of the interface to the edges of the mesh}. This is also illuminated in Figure \ref{fig:straightangled0} where we see that the angle between $S_{FEM}^h$ and $S_{ENR}^h$ for the SGFEM approaches a fixed value bounded away from 0. However, similarly to the 1-D interface problem, there has to be a safety-check and it is advisable not to enrich a node if the interface is very close to it; otherwise, round-off errors could contaminate the solution.

We summarize the results mentioned above about different GFEMs in the following table.

\begin{center}
\begin{tabular}{c|ccccc}
\ & FEM & GFEM (Topo) & GFEM (Geo) & GFEM (M) & SGFEM \\ \hline
Order of Conv. & $O(h^{1/2})$ & $O(h^{1/2})$ & $O(h)$ & $O(h)$ & $O(h)$ \\ \hline
\multirow{2}{*}{Angle} & & bounded away & $\to 0$ & $\to 0$ only & bounded away \\ 
& & from 0 & as $h \to 0$ & as $\Gamma \to$ edge & from 0 \\ \hline
$\kappa_2(\bm{A})$ & $O(h^{-2})$ & $O(h^{-2})$ & $O(h^{-4})$ & $O(h^{-2})$ & $O(h^{-2})$ \\ \hline
Robustness & yes & yes & yes & no & yes
\end{tabular}
\end{center}

We thus conclude that among all the GFEMs considered in this section, the SGFEM is the only method that has all the desired features -- it yields accurate approximation, it is well-conditioned, and it is robust.

\section{Circular interface problem}\label{sec5}

In this section, we discuss the GFEM applied to a 2-D problem with a circular interface. We consider a specific manufactured problem such that the solution does not have any singularities. This problem will be solved by extending what was done in Section 4 on a straight interface to a circular interface.

\begin{figure}[!ht]
\centerline{\includegraphics[scale = 0.65]{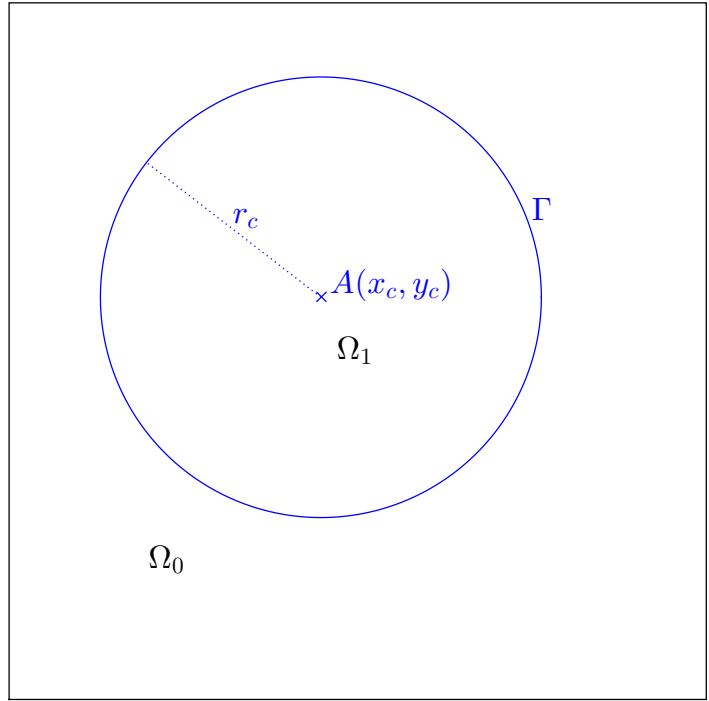}}
\caption{Circular interface problem.}
\label{fig:circlepblm}
\end{figure}

Consider the domain $\Omega = (0,1) \times (0,1)$. For given $r_c>0$ and $(x_c,y_c) \in \Omega$, let $\Gamma:= \{\bm{x} \in \Omega\,:\, \gamma(\bm{x})=0\}$ be the interface, where $\gamma(\bm{x})=(x-x_c)^2 + (y-y_c)^2 - r_c^2$ is the circle of center $A(x_c,y_c)$ and radius $r_c$, as shown in Figure \ref{fig:circlepblm}, where we have chosen $r_c = 1/\sqrt{10}$ and $(x_c,y_c) = (1/\sqrt{5},1/\sqrt{3})$. We set $\Omega_0:= \Omega \cap \{\bm{x}\,:\, \gamma(\bm{x}) <0\}$ and $\Omega_1:= \Omega \cap \{\bm{x}\,:\, \gamma(\bm{x}) >0\}$. Note that when varying $x_c,y_c,r_c$, the interface $\Gamma$ varies as well.

We consider the problem \eqref{WFneu} with $f\equiv 0$, $q \equiv 0$, and $g_N$ satisfying the compatibility condition \eqref{CompCond}. Moreover, we consider $a_i(\bm{x}) = a_i$, $i=0,1$, where $a_0, a_1$ are strictly positive constants. The solution $u \in \mc{E} = H^1(\Omega)$ of \eqref{WFneu} exists and is unique up to an additive constant. We set $\|u\|_{\mc{E}}:= B(u,u)^{1/2}$.

For the computations presented in this section, we will consider the manufactured solution of \eqref{WFneu}, given by
\begin{equation}
u_{ex} = \left\{ \begin{array}{lr}
r^2 \cos(2\theta) +C,& r \le r_c,\\
B_0 r^2 \cos(2\theta) + B_1 r^{-2} \cos(2\theta) +C,& r \ge r_c,\end{array} \right. \label{ManufacSol2}
\end{equation}
where $(r,\theta)$ is the polar coordinate centered at $A$. We choose $B_0,B_1$ such that $u_{ex}$ and $a(\bm{x})\, \frac{\partial u_{ex}}{\partial n}$ are continuous across the interface $\Gamma\, (r = r_c)$ and then $C$ so that $u_{ex}(0,0) = 0$. It is clear that $u_{ex}$ is continuous in $\Omega$ with no singularity in $\ov{\Omega}$. Moreover $f = 0$ and $g_N$, in \eqref{WFneu}, is obtained from $u_{ex}$ using \eqref{BCneu}.

\begin{figure}[!ht]
\centerline{\includegraphics[scale = 0.65]{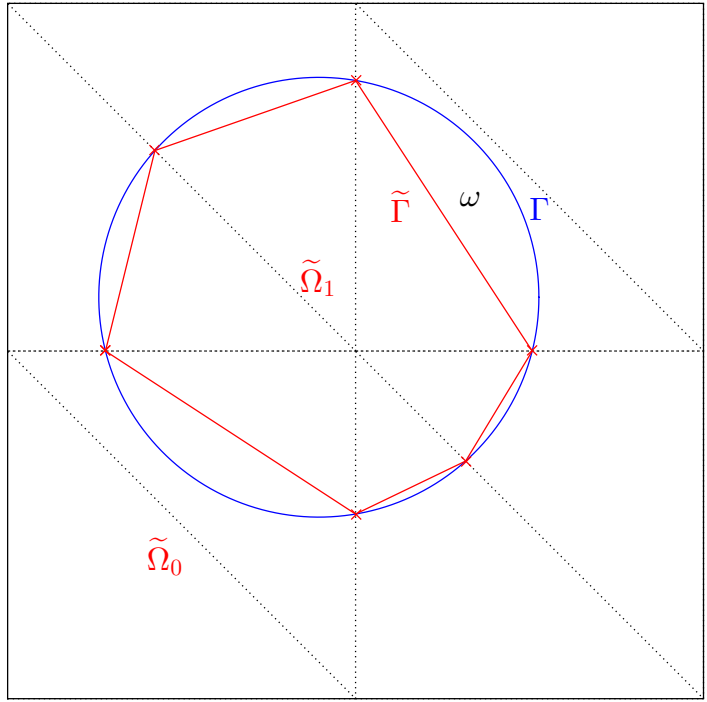}}
\caption{Circular interface and perturbed interface.}
\label{fig:circleperturb}
\end{figure}

Let us first comment on FEM. We use the same discretization and notations as in Section \ref{sec4} to define $\bm{x_i}, \tau, \omega_{\bi}, N_{\bi}$ with an understanding that they depend on $h$. During the assembling of the stiffness matrix $\bm{A}_{11}$ we need to evaluate quantities of the form
\begin{equation*}
B(N_{\bm{j}},N_{\bi}) = \int_{\Omega} a \nabla N_{\bm{j}} \nabla N_{\bi} \, d\bm{x}.
\end{equation*}
With $a$ being discontinuous on $\Gamma$, we thus have to split the integration domain into two complementary parts, one in $\Omega_0$ and one in $\Omega_1$. Each of these sub-parts is then no longer polygonal as the interface is curved, numerical integration up to machine precision would thus be costly and difficult to implement. We propose another approach instead. Considering that the interface can be described by $\gamma(\bm{x})=(x-x_c)^2 + (y-y_c)^2 - r_c^2$, we can readily invert this expression to find the intersections of the finite element triangulation with the interface. This gives a set of points whose convex hull forms a polygon, noted $\widetilde{\Gamma}$ -- see Figure \ref{fig:circleperturb}, where $h = 1/2$, $r_c = 1/\sqrt{10}$ and $(x_c,y_c) = (1/\sqrt{5},1/\sqrt{3})$. We will use this polygon instead of the circular interface.

We define a ``perturbed'' $a(\bm{x})$, noted $\widetilde{a}(\bm{x})$, whose value is $a_0$ outside the polygon, domain denoted $\widetilde{\Omega}_0$ and $a_1$ inside, domain denoted $\widetilde{\Omega}_1$. We also define $\omega$ as the difference $\omega :=  \widetilde{\Omega}_0 \backslash \Omega_0$ (we could equivalently define $\omega$ as the difference $\Omega_1 \backslash \widetilde{\Omega}_1$). This approach could be considered as a perturbation of the original problem. We can now see why we first studied the straight interface problem. With the perturbed interface, instead of \eqref{InterfaceProb1d} we now consider the perturbed variational problem: find $\widetilde{u} \in \mc E$ satisfying
\begin{align}
\widetilde{B}(\widetilde{u},v) = F(v), \quad \mbox{for all }v \in \mc E, \label{Interfaceperturb}
\end{align}
where
\begin{align*}
\widetilde{B}(u,v) := \int_\Omega \widetilde{a} \nabla u \cdot \nabla v\, d \bm x. 
\end{align*}
Note that while the solution $u \in \mc E$ \eqref{ManufacSol2} of the original problem \eqref{InterfaceProb1d} does not have any singularity in $\ov \Omega$, the perturbed solution $\widetilde{u} \in \mc E$ of the perturbed problem \eqref{Interfaceperturb} may exhibit singularities due to the corners of the perturbed interface $\widetilde{\Gamma}$. However, we do not solve \eqref{Interfaceperturb}. Instead, we solve the finite dimensional problem: find $u_h \in S^h$ satisfying
\begin{align}
\widetilde{B}(u_h,v) = F(v), \quad \mbox{for all }v \in S^h, \label{Interfaceperturbfinitedim}
\end{align}
so the solution $u_h \in S^h$ does not have any singularity in $\ov \Omega$.


Since numerical solutions $u_h \in S^h$ of \eqref{Interfaceperturbfinitedim} satisfy Galerkin orthogonality only with respect to $\widetilde{u}$ and $\widetilde{B}(\cdot,\cdot)$ and not with respect to $u$ and $B(\cdot,\cdot)$, we compute the discretization error as $\epsilon^h:=\left| \|u\|^2_\mc{E}-\|u_h\|^2_\mc{\widetilde{E}}\right|^{1/2}$, where $\|v\|_\mc{\widetilde{E}} := \widetilde{B}(v,v)^{1/2}$. First, the energy of the exact solution $u \in \mc E$ is still computed with respect to $B(\cdot,\cdot)$ (i.e., the true interface $\Gamma$), while the energy of the discrete solution $u_h \in S^h$ is now computed with respect to $\widetilde{B}(\cdot,\cdot)$ (i.e., the perturbed interface $\widetilde{\Gamma}$). Secondly, we compute the ``difference of the energy norms'' and not ``energy norm of the difference'' (note that the two usually coincide thanks to Galerkin orthogonality). We use this unusual definition of the discretization error for computational reasons. This definition avoids integrating the quantity $a \nabla (u-u_h) \cdot \nabla (u-u_h)$ on each element $\tau \in E$, which would be costly and difficult to implement due to the curved interface $\Gamma$. Instead, we thus compute the discretization error as $\epsilon^h=\left| \|u\|^2_\mc{E}-\|u_h\|^2_\mc{\widetilde{E}}\right|^{1/2}$, and it holds
\begin{align*}
\left| \| u -u_h \|^2_{\mc E} - \left|\|u\|^2_\mc{E}-\|u_h\|^2_\mc{\widetilde{E}}\right| \right| &\leq O(h^k),
\end{align*}
where $k = 3/2$ for FEM and $k = 2$ for GFEM \& SGFEM. The proof of this result can be found in Appendix \ref{app:error}. Note that we have committed two crimes: first, the perturbation of the interface from $\Gamma$ to $\widetilde{\Gamma}$, secondly, the computation of the error as $\left| \|u\|^2_\mc{E}-\|u_h\|^2_\mc{\widetilde{E}}\right|^{1/2}$ and not as $\| u -u_h \|_{\mc E}$. However, these two crimes have limited effects compared to the true discretization error $\| u -u_h \|_{\mc E}$ which is of order $O(h^{1/2})$ for FEM and $O(h)$ for GFEM \& SGFEM.

Finally, note that the quantity $\|u\|^2_\mc{E}$ in $\epsilon^h$ can be calculated as $\int_{\partial \Omega} u g_N\, ds$, thus avoiding the curved interface $\Gamma$.

Let us now describe the GFEM in 2-D. The approximation space $S^h$ of the GFEM is once again given by $S^h = S_{FEM}^h \oplus S_{ENR}^h$.

\begin{figure}[!ht]
\centerline{\includegraphics[scale = 0.65]{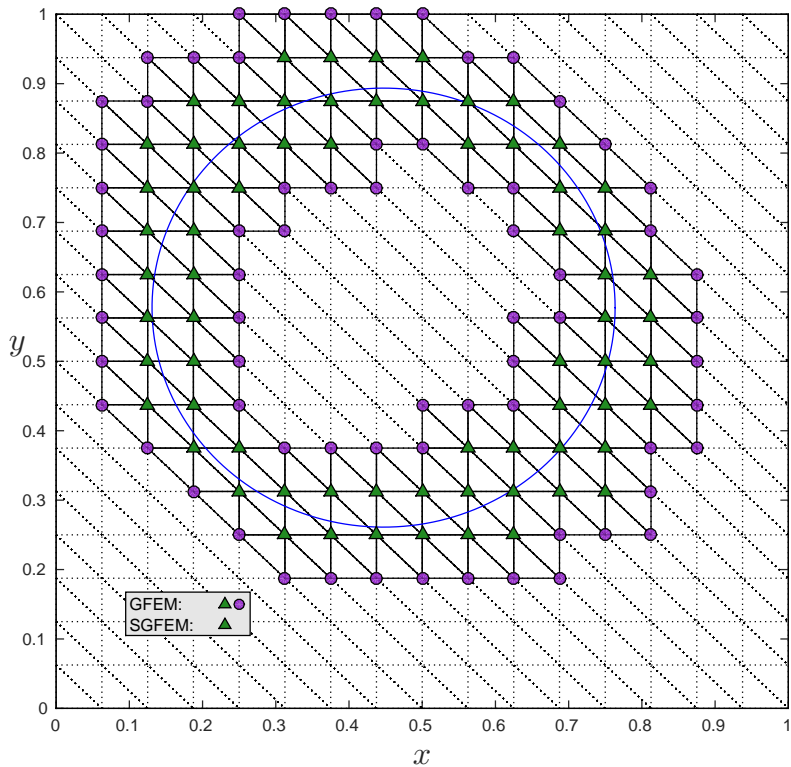}}
\caption{Nodes enriched for the circular interface problem.}
\label{fig:circlepatches}
\end{figure}

In this section, we will only consider the so-called M-GFEM of the previous section. Recall that Topological GFEM does not recover the optimal rate of convergence in terms of error, while Geometric GFEM is badly conditioned; we thus do not discuss them anymore. \textit{\ul{For brevity, we will refer to M-GFEM simply as GFEM in this section and the next}}. The enrichment function $w$ for the circular problem is again based on the so-called \textit{distance function}. However, note that $\mbox{dist}(\bm{x},\Gamma)$ is quadratic in $\Omega_0$ and $\Omega_1$. Moreover, $\mbox{dist}(\bm{x},\widetilde{\Gamma})$ is quadratic in some regions of $\Omega$. In order to facilitate numerical integration, we would like to use a piecewise linear enrichment function. Here is how we were able to obtain such a function. We start with $w^\diamond(\bm{x}) = \mbox{dist}(\bm{x},\Gamma)$ the distance to the interface $\Gamma$. Next, we perform a triangulation of $\mc{T}_h$ using $\widetilde{\Gamma}$ as an edge constraint: each element in $E_\Gamma$ is divided into elementary triangles whose edges do not cross $\widetilde{\Gamma}$. We then compute the linear interpolant of $w^\diamond(\bm{x})$ on this triangulation, thus giving $w^*(\bm{x})$. It is continuous and piecewise linear in $\Omega$. We mention the presence of ``shadow interfaces'' due to the triangulation, but the additional weak discontinuities (kinks) are controlled by the true distance to the interface rather than by arbitrary factors. Then, the enrichment function used in GFEM is given by
\begin{eqnarray*}
&&w(\bm{x}) = \left\{ \begin{array}{ll}
w^*(\bm{x}), & \bm{x} \in \tau \in E_\Gamma, \\
\mbox{linear function}, & \bm{x} \in \tau \in E \backslash E_\Gamma,
\end{array} \right. \\
&& w(\bm{x_i}) = 0, \quad \bm{x_i} \mbox{ is not a vertex of } \tau \in E_\Gamma.
\end{eqnarray*}
Furthermore, the set of indices of enriched nodes is given by
\[
\mc{R}^h=\{\bi \in \mc{N}^h\,:\, \bm{x_i} \in \bigcup_{\omega_{\bm{j}} \cap \Gamma \ne \emptyset} \ov{\omega}_{\bm{j}} \}.
\]
The enriched nodes $\bm{x_i},\bi \in \mc{R}^h$ are shown in Figure \ref{fig:circlepatches}, where $h=1/16$, $r_c = 1/\sqrt{10}$ and $(x_c,y_c) = (1/\sqrt{5},1/\sqrt{3})$. Note that like in Section \ref{sec4}, card$\{\bm{x_i}\}_{\bi \in \mc{R}^h} = O(h^{-1})$.

\begin{figure}[!ht]
\centerline{\includegraphics[scale = 0.65]{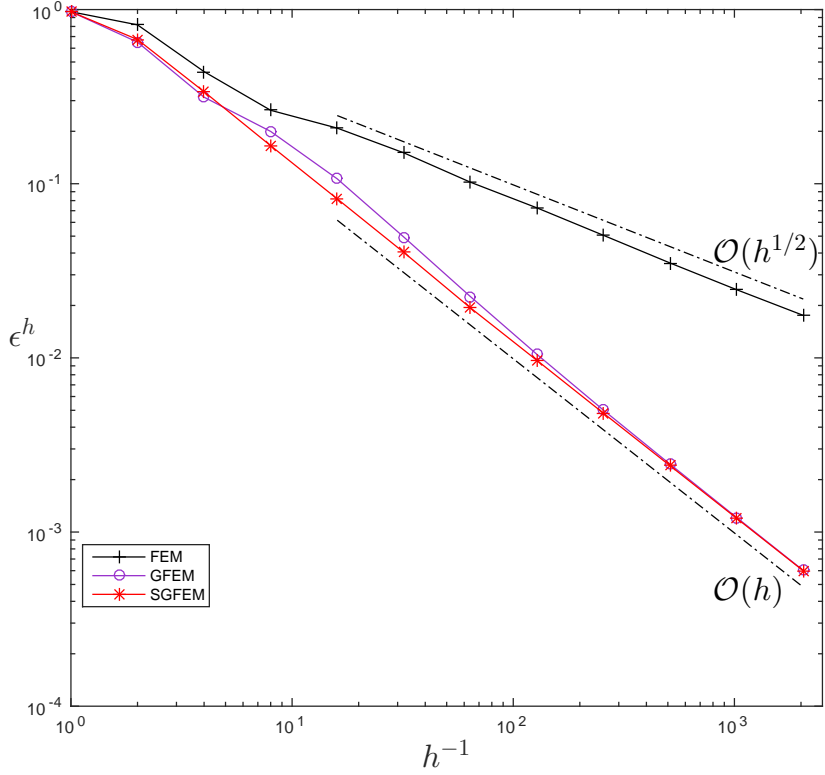}}
\caption{Relative error in the energy norm against $h$.}
\label{fig:circleerrorh}
\end{figure}

First, let us consider the exact solution $u\in \mc{E}$ of \eqref{WFneu} given in \eqref{ManufacSol2} with $r_c=1/\sqrt{10}$, $(x_c,y_c) = (1/\sqrt{5},1/\sqrt{3})$, $a_0 = 1$ and $a_1 = 10$. We computed the error $\epsilon^h = \left|\|u\|^2_\mc{E}-\|u_h\|^2_\mc{\widetilde{E}}\right|^{1/2}$, where $u_h$ is the solution of \eqref{Interfaceperturbfinitedim} associated with the GFEM described above. The log-log plot of the (relative) error is given in Figure \ref{fig:circleerrorh}. It is clear that GFEM yields the convergence of $O(h)$, whereas the order of convergence for standard FEM is only $O(h^{1/2})$. This is very similar to the results of Section \ref{sec4}.

\begin{figure}[!ht]
\centerline{\includegraphics[scale = 0.65]{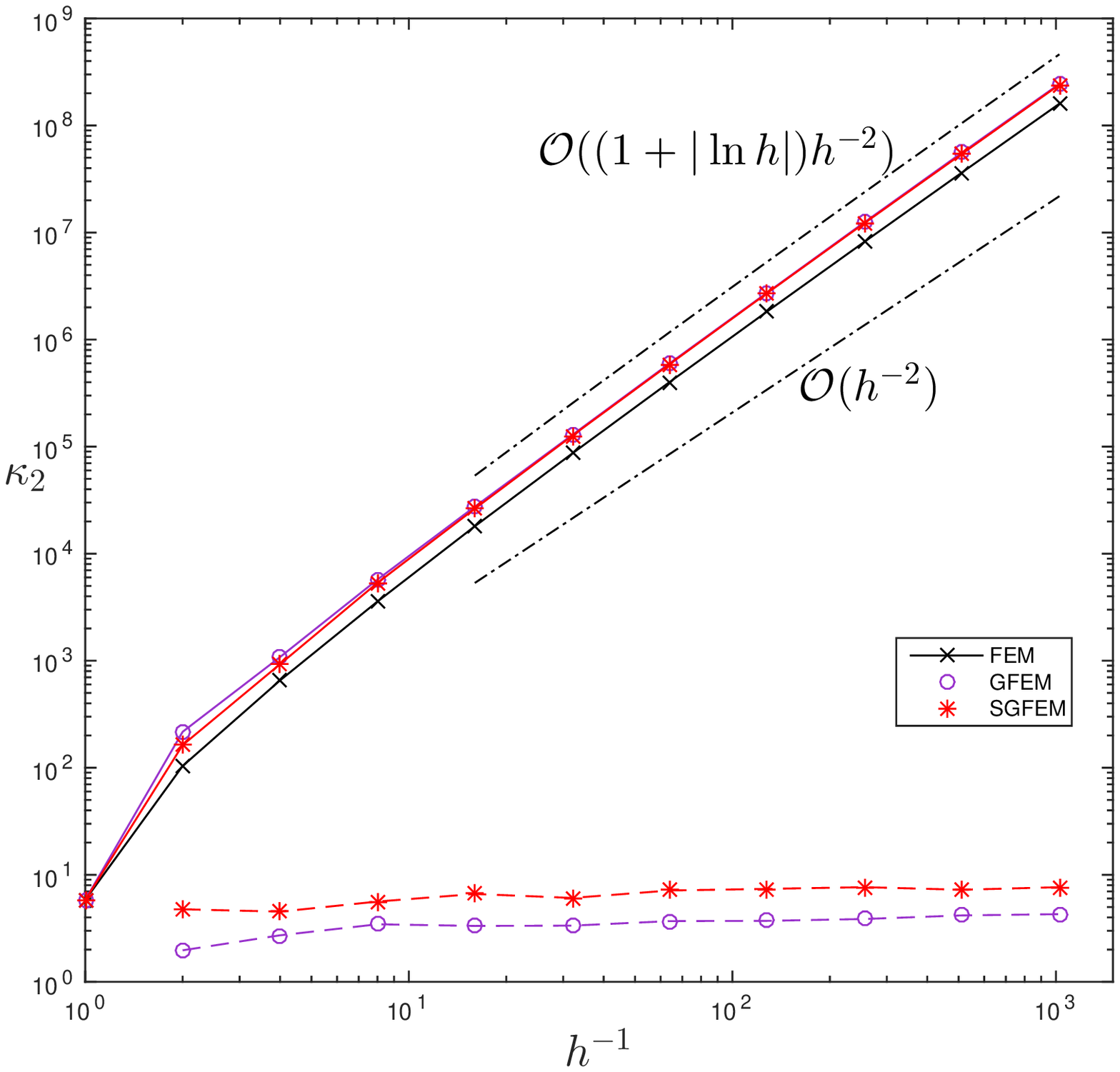}}
\caption{Condition number of the scaled stiffness matrices $\bm{A}$ (solid lines) and $\bm{A}_{22}$ (dashed lines) against $h$.}
\label{fig:circlecondh}
\end{figure}

Concerning the condition number of the scaled stiffness matrix $\bm{A}$ in \eqref{goodLinSys} for different values of $h$, we display the results in Figure \ref{fig:circlecondh} (solid lines), again with $r_c=1/\sqrt{10}$, $(x_c,y_c) = (1/\sqrt{5},1/\sqrt{3})$, $a_0 = 1$ and $a_1 = 10$. It is clear that for GFEM, the condition number $\kappa_2(\bm{A}) = O(h^{-2})$. We mention that $\kappa_2(\bm{A}_{11}) = O(h^{-2})$, where $\bm{A}_{11}$ is the stiffness matrix of the standard FEM. In other words, the conditioning of GFEM is of the same order as that of a standard FEM. Again, this is very similar to the results in Section \ref{sec4}. We also mention that $\kappa_2(\bm{A}_{22})$ is bounded for all values of $h$ (see dashed lines in Figure \ref{fig:circlecondh}).

\begin{figure}[!ht]
\centerline{\includegraphics[scale = 0.65]{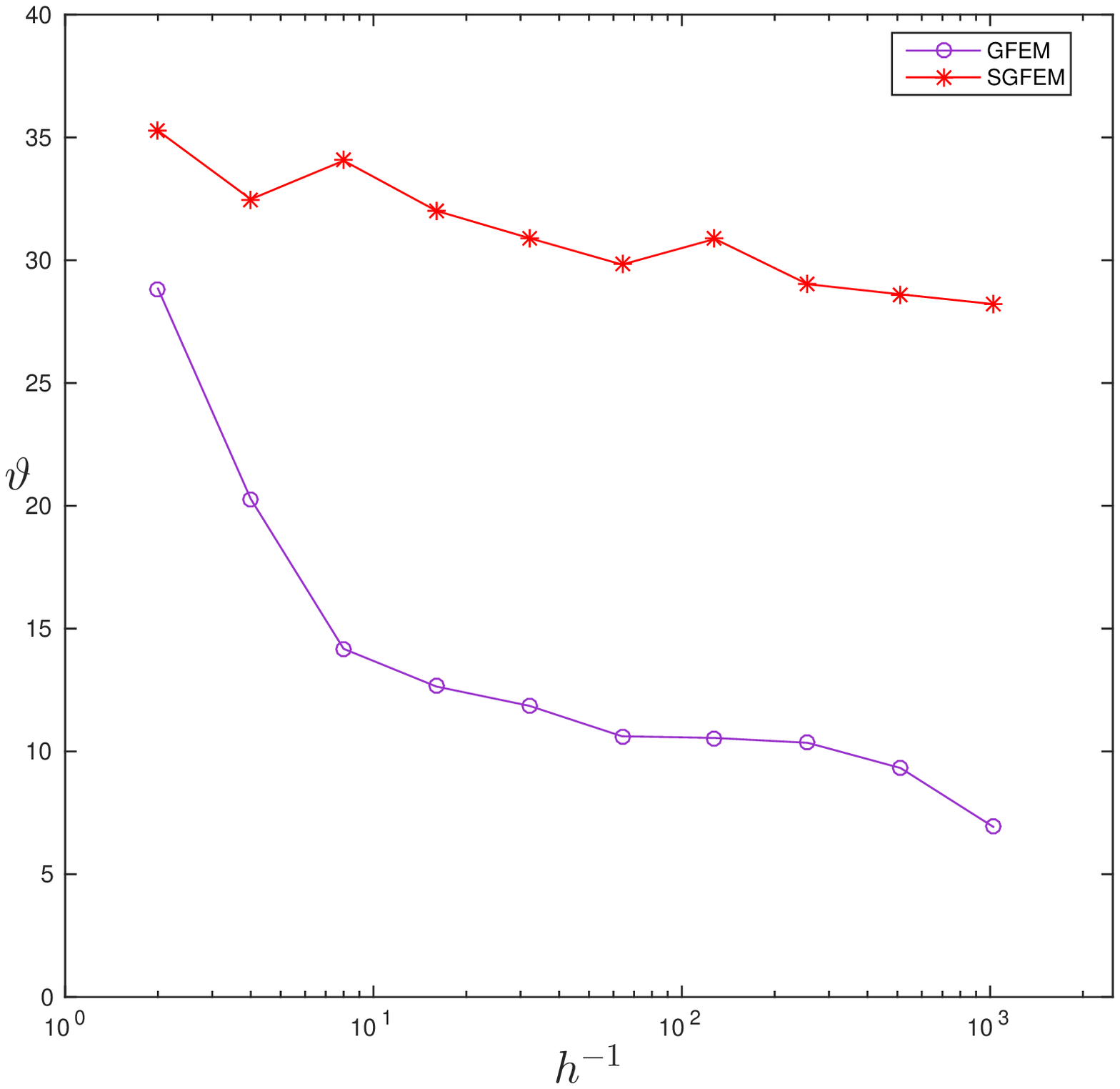}}
\caption{Angle (in degrees) between the spaces $S_{FEM}^h$ and $S_{ENR}^h$ against $h$.}
\label{fig:circleangleh}
\end{figure}

We have computed the angle for the GFEM for different values of $h$ and displayed the results in Figure \ref{fig:circleangleh}, again with $r_c=1/\sqrt{10}$, $(x_c,y_c) = (1/\sqrt{5},1/\sqrt{3})$, $a_0 = 1$ and $a_1 = 10$. It is clear that the angle remains bounded away from $0$ for all the values of $h$ and, it thus illuminates the result presented above since we have seen in Figure \ref{fig:circlecondh} that $\kappa_2(\bm{A})= O(h^{-2})$.

\textbf{Stable GFEM (SGFEM):} The enrichment is defined as
\begin{equation}
w(\bm{x}) = w^*(\bm{x}) - \mc{I}_hw^*(\bm{x}), \label{SGFEMenr2}
\end{equation}
where $\mc{I}_hw^*(\bm{x})$ is the piecewise linear interpolant of $w^*(\bm{x})$ with respect to the triangulation $\mc T_h$. Since $w^*(\bm{x})$ is linear on $\tau \in E \backslash E_\Gamma$, it is easy to see that 
\[
\mbox{supp}\{w(\bm{x})\} = \bigcup_{\tau \in E_\Gamma} \tau.
\]
Moreover, the indices of the enriched nodes are 
\begin{equation}
\mc{R}^h = \{\bi\,:\, \omega_{\bi} \cap \Gamma \ne \emptyset \}. \label{SGFEMenrNodes2}
\end{equation}
In Figures \ref{fig:circleerrorh}, \ref{fig:circlecondh} and \ref{fig:circleangleh}, we have plotted the error $\epsilon^h=\left| \|u\|^2_\mc{E}-\|u_h\|^2_\mc{\widetilde{E}}\right|^{1/2}$, the condition number $\kappa_2(\bm{A})$, and the angle between $S_{FEM}^h$ and $S_{ENR}^h$, with respect to $h$, for the GFEM with enrichment given in \eqref{SGFEMenr2} and $\mc{R}^h$ as in \eqref{SGFEMenrNodes2}. It is clear that the method yields the optimal order of convergence of $O(h)$, $\kappa_2(\bm{A}) = O(h^{-2})$, and the angle between $S_{FEM}^h$ and $S_{ENR}^h$ is bounded away from $0$ for all the values of $h$ considered in the experiment. Condition \eqref{AngleCond1} is thus satisfied. We have also checked that \eqref{Assump2} is satisfied (see dashed lines in Figure \ref{fig:circlecondh}). The GFEM with enrichment given in \eqref{SGFEMenr2}--\eqref{SGFEMenrNodes2} is thus an SGFEM. Once again, the angle between $S_{FEM}^h$ and $S_{ENR}^h$ for the SGFEM is larger than that for the GFEM. This feature is central to the iterative solution of \eqref{goodLinSys}, which we will illustrate in Section \ref{sec6}. We thus conclude that the SGFEM once again has all the desired features: it yields accurate approximation, it is well-conditioned, and it is robust.

We mention that if the closed curved interface $\Gamma$ has a straight part, the angle between the spaces $S_{FEM}^h$ and $S_{ENR}^h$ associated with the GFEM considered in this section will become small when the relative distance between the element edges and the ``straight part of $\Gamma$" is small, similar to what we observed in Figure \ref{fig:straightangled0}. This phenomenon will give rise to a much larger value of $\kappa_2(\bm{A})$ similar to the situation shown in Figure \ref{fig:straightcondd0}. The GFEM will thus not be stable. The SGFEM will nevertheless be stable for such interface problems. However, since in this section we considered a circular interface problem, which does not have any straight parts, the GFEM did not exhibit a behavior similar to what was observed in Figure \ref{fig:straightcondd0} or Figure \ref{fig:straightangled0}.

\section{Iterative methods}\label{sec6}
In this section, we exploit the angle condition between $S_{FEM}^h$ and $S_{ENR}^h$ and design appropriate iterative solvers based on previous observations. For a given error tolerance, we will compare the performance of the iterative solvers for FEM, GFEM (recall that we are only referring to M-GFEM; see Section \ref{sec5}) and SGFEM. The solvers designed in this section can be applied to both the straight interface problem (Section \ref{sec4}) and the circular interface problem (Section \ref{sec5}) and the conclusions we come to are very similar for the two problems. Note that the discretization error takes a different form depending on the problem: for the straight interface case, it is the classical $\epsilon^h = \|u - u_h\|_{\mc{E}}$, while for the circular interface case we consider instead $\epsilon^h = \left|\|u\|^2_{\mc{E}}-\|u_h\|^2_{\mc{\widetilde{E}}}\right|^{1/2}$. Let $v_h \in S^h$ be an approximation of $u_h \in S^h$. \textit{For the sake of concision, we will refer to the truncation error as $\delta := \|u_h - v_h\|_{\mc{E}}$} with the understanding that in the case of the circular interface problem, we actually mean $\delta = \|u_h - v_h\|_{\mc{\widetilde{E}}}$.

We start by examining the linear system associated to standard FEM \eqref{goodLinSys} with $w \equiv 0$, yielding $S^h = S_{FEM}^h$ and $ \bm{A}= \bm{A}_{11}$. The exact (discrete) solution of this system is noted $u_h$. We design an iterative solver based on Schur complement \cite{Briggs2000, Hackbusch1985}. We will denote by $v_h^{i}$ the iterative solution at iteration $i$. There are now two sources of error: the first one is the discretization error $\epsilon^h$ due to the choice of the approximation space $S^h$. The second one is the truncation error due to the choice of the iterative solver, $\delta_i = \|u_h - v_h^{i}\|_{\mc{E}}$ (or $\delta_i = \|u_h - v_h^{i}\|_{\mc{\widetilde{E}}}$ for the circular interface problem). The iterative solver we will use is Conjugate Gradient (CG) preconditioned by Full Multigrid (FMG). First, we recall the essential steps in multigrid methods, when they are viewed not as preconditioner but as solvers.

We start by applying a few relaxation steps (e.g., Gauss-Seidel) until the residual starts to stagnate, indicating that the low frequency content of the solution has been found. Then, we interpolate the residual onto a coarser grid and perform again a few relaxations until the residual stagnates again. This is applied until the coarsest level is reached (typically containing only a few elements), where we solve the residual equation exactly, before projecting back onto the finer grids, applying there again a few relaxations each time. This scheme is the so-called V-cycle. Applying successive V-cycles allows the truncation error $\delta_i$ to decrease geometrically from one cycle to the next. FMG is a variant of this scheme. In FMG, V-cycles are applied recursively starting on the coarsest mesh, until the finest is reached. Once again, when several FMGs are performed, the truncation error $\delta_i$ decreases geometrically, with a higher reduction factor than for a single V-cycle. However, the computational cost of an FMG is slightly larger than that of a V-cycle \cite{Briggs2000}. 

As stated before, we will not use FMG as a solver in this paper, but rather as a preconditioner for CG. One of the reasons for this choice is that the mathematical theory for multigrid methods is available for smooth coefficients, but it is not yet well-developed for discontinuous coefficients. As a result, we have observed that CG preconditioned by a multigrid method such as V-cycle or FMG was more robust than the multigrid method alone used as a solver. We have also noted that FMG has better computational properties than the regular V-cycle (reduction of the truncation error for the same computational work). As a result, our solver for $S_{FEM}^h$ is CG preconditioned by FMG.

We now have to design an appropriate stopping criterion in order to decide when to stop the CG iterations. By a priori error estimation, we know that the discretization error $\epsilon^h$ behaves in some $Ch^p$, where $p$ is known by the underlying properties of the PDE, the choice of the partition of unity and the choice of the enrichment. A wise stopping criterion would be to stop the iterations as soon as the truncation error $\delta_i$ becomes significantly smaller than the discretization error $\epsilon^h$. We will show in Appendix \ref{app:stoppingcrit} how to attain this by using a controlling factor. Performing more iterations will not result in a sensible improvement as the quality of the iterative solution $v_h^i$ will mostly be driven by the discretization and not the truncation part of the error.

Due to the use of the CG solver, the truncation error is no longer expected to decay geometrically with the number of iterations. However, the use of FMG as a preconditioner reduces the ``effective" condition number $\kappa_2 (\bm{A}_{11})$ from $O(h^{-2})$ to $O(h^{-1})$, see \cite{Johnson2012}. We thus rely on an error estimator for the truncation error $\delta_i$ based on the residual of the linear system and the ``effective" spectral radius of the matrix $\bm{A}_{11}^{-1}$ from \eqref{goodLinSys} using a so-called inverse estimate. This allows us to estimate the truncation error at step $i$ as follows
\begin{align}
e^i = \frac{\|\bm{f}_1 -\bm{A}_{11}\bm{x}^i\|_{l^2}}{h}, \label{eq:estimfem}
\end{align}
where the vector $\bm{x}^i = \bm{D}^{-1} \hat{\bm{c}}^i$ is associated to $v_h^i$ in the approximation space $S^h_{FEM}$, i.e., $v_h^i = \sum_{k\in \mc{N}^h_d} c_k^iN_k$, and $\bm{D}$ is the scaling matrix associated to $\bm{A}$ from \eqref{goodLinSys}. The full derivation of this estimator can be found in Appendix \ref{app:stoppingcrit}.

The algorithm developed for FEM schematically takes the form of Algorithm \ref{alg:fem}.

\begin{algorithm}[!ht]
 \KwData{$h,\bm{A}_{11},\bm{f}_1,k$}
 \KwResult{$v^{i^*}_h,i^*$}
 $\epsilon = h^{1/2}, v^0_h=0, e^0=\infty, i=0$\;
 \While{$e^i \geq \epsilon/k$}{
 	 $i \leftarrow i+1$\;
     Compute $v^i_h$ using initialization $v^{i-1}_h$\;
     Compute error estimator $e^i$ using \eqref{eq:estimfem}\;
 }
 $i^* = i$.
  \caption{Algorithm for FEM.}
 \label{alg:fem}
\end{algorithm}

Note that, as always, we work on the scaled system \eqref{goodLinSys}. The algorithm stops iterating as soon as the estimated truncation error $e^i$ becomes significantly lower than the a priori estimated discretization error $\epsilon$. The factor $k$ is used to control the different constants of proportionality appearing in the intermediate calculations (see Appendix \ref{app:stoppingcrit}).

\begin{figure}[!ht]
\centerline{\includegraphics[scale = 0.5]{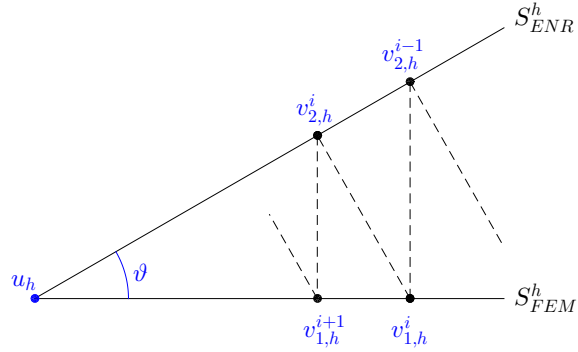}}
\caption{Outside iterations: block Gauss-Seidel scheme between $S^h_{FEM}$ and $S^h_{ENR}$ \eqref{eq:Sfem}--\eqref{eq:Senr}.}
\label{fig:gaussseidel0}
\end{figure}

Let us now consider the case of GFEM \& SGFEM. Since we want to exploit the angle condition between $S_{FEM}^h$ and $S_{ENR}^h$, the idea is to apply a block Gauss-Seidel iterative scheme between $S_{FEM}^h$ and $S_{ENR}^h$. This defines our ``outer iterations''. The system to be solved is \eqref{goodLinSys}, and by partitioning the solution $\bm{x}$ in $\bm{x} = \left[ \bm{x}_1, \bm{x}_2 \right] = \left[ x_{1,i}, x_{2,k} \right]_{i\in \mc{N}^h_d;k\in \mc{R}^h}$, corresponding to $u_h \in S^h = S_{FEM}^h \oplus S_{ENR}^h$ where $u_h = u_{1,h} + u_{2,h}$ with $u_{1,h} \in S_{FEM}^h$ and $u_{2,h} \in S_{ENR}^h$, we obtain
\begin{eqnarray}
&\mbox{Solve }\bm{A}_{11} \bm{x}^i_1 = \bm{f}_1 - \bm{A}_{12} \bm{x}^{i-1}_2 \mbox{ for } \bm{x}^i_1,\label{eq:Sfem}\\[5pt]
&\mbox{Solve }\bm{A}_{22} \bm{x}^i_2 = \bm{f}_2 - \bm{A}_{12}^T \bm{x}^i_1 \mbox{ for } \bm{x}^i_2,\label{eq:Senr}
\end{eqnarray}
where $\bm{x}^i_\cdot$ denote successive iterates $\bm{x}^i = \left[ \bm{x}^i_1, \bm{x}^i_2 \right]$ corresponding to iterates in $S^h$, $v_h^i = v^i_{1,h} + v^i_{2,h}$ with $v^{i}_{1,h} \in S_{FEM}^h$ and $v^{i}_{2,h} \in S_{ENR}^h$. The truncation error $\delta_i = \|u_h - v_h^{i}\|_{\mc{E}}$ decreases geometrically with a ratio related to the angle between the spaces $S_{FEM}^h$ and $S_{ENR}^h$. In fact, if $q$ denotes this ratio, we have $q =\cos^2\left(\vartheta(S^h_{FEM},S^h_{ENR})\right)$, see Figure \ref{fig:gaussseidel0} where $\vartheta = \pi/6$ and the truncation error is divided at each iteration by a factor $q^{-1}=4/3$, as are the quantities $\| v_{1,h}^i -v_{1,h}^{i-1} \|_{\mc{E}}$, $\| v_{2,h}^i -v_{2,h}^{i-1} \|_{\mc{E}}$ and $\| v_h^i -v_h^{i-1} \|_{\mc{E}}$.

We can estimate the truncation error in the outer iteration $\delta_i = \|u_h - v_h^i\|_{\mc{E}}$ using Richardson extrapolation on the last three iterates as follows
\begin{align}
e^i = \frac{1}{\displaystyle\frac{1}{\|v_h^i - v_h^{i-1}\|_{\mc{E}}}-\frac{1}{\|v_h^{i-1} - v_h^{i-2}\|_{\mc{E}}}}. \label{eq:gscheat}
\end{align}
The derivation of this estimator can be found in Appendix \ref{app:stoppingcrit}.

At a given stage $i$ in the outer iteration, we now have to find $v^{i}_{1,h}$ using \eqref{eq:Sfem}, and $v^{i}_{2,h}$ using \eqref{eq:Senr}. Each of these will also be done in an iterative manner, which defines our ``inner iterations".

For $v^{i}_{1,h}$, i.e., \eqref{eq:Sfem}, we have to invert $\bm{A}_{11}$, we thus use the same solver as in FEM, that is CG preconditioned by FMG, building a sequence $w^{i,j}_{1,h}$ governed by iteration number $j$. The stopping criteria need not be as demanding as in FEM because we are not interested in $v^{i}_{1,h}$ but in $u_{1,h}$ instead. We thus design the stopping criterion so that the truncation error in the inner iteration in $S_{FEM}^h$, $\delta^{i,j}_1 = \|v^{i}_{1,h} - w^{i,j}_{1,h}\|_{\mc{E}}$ is only a small fraction of the truncation error in $S^h$, $\delta_i = \|u_h - v^{i}_h\|_{\mc{E}}$. This yields the following error estimator for $\delta^{i,j}_1 = \|v^{i}_{1,h} - w^{i,j}_{1,h}\|_{\mc{E}}$
\begin{align}
e^{i,j}_1 = \frac{\|\bm{f}_1-\bm{A}_{12}\bm{x}_2^{i-1}-\bm{A}_{11}\bm{x}_1^{i,j}\|_{l^2}}{h}. \label{eq:innerSfem}
\end{align}
The derivation of this estimator essentially follows the same as for \eqref{eq:estimfem}, having replaced the right hand side $\bm{f}_1$ with $\bm{f}_1-\bm{A}_{12}\bm{x}_2^{i-1}$. Once the stopping criterion is reached, after say $j'$ iterations, we set $w_{1,h}^i := w_{1,h}^{i,j'}$ and update $w_h^i = w_{1,h}^i + w_{2,h}^i$ (in practice $w_{2,h}^i$ is to be computed in the second block of the Gauss-Seidel scheme). However, $v_h^i$ is not available in practice, $\delta_i$ is therefore not estimated using \eqref{eq:gscheat} but rather
\begin{align}
e^i = \frac{1}{\displaystyle\frac{1}{\|w_h^i - w_h^{i-1}\|_{\mc{E}}}-\frac{1}{\|w_h^{i-1} - w_h^{i-2}\|_{\mc{E}}}}, \label{eq:gsnocheat}
\end{align}
which is the same as \eqref{eq:gscheat}, having replaced $v_h^\cdot$ with $w_h^\cdot$.

For $v^{i}_{2,h}$, i.e., \eqref{eq:Senr} we have to invert $\bm{A}_{22}$, and based on condition \eqref{Assump2}, its condition number is bounded, so we choose to solve the equation in $S_{ENR}^h$ using CG. We denote the successive iterates $w^{i,j}_{2,h}$. Once again, we design the stopping criterion so that the truncation error in the inner iteration in $S_{ENR}^h$, $\delta^{i,j}_2 = \|v^{i}_{2,h} - w^{i,j}_{2,h}\|_{\mc{E}}$ is only a small fraction of the truncation error in $S^h$, $\delta_i = \|u_h - v^{i}_h\|_{\mc{E}}$. To estimate $\|v^{i}_{2,h} - w^{i,j}_{2,h}\|_{\mc{E}}$, we use the norm of the residual. This yields the following error estimator for $\delta^{i,j}_2 = \|v^{i}_{2,h} - w^{i,j}_{2,h}\|_{\mc{E}}$
\begin{align}
e^{i,j}_2 = \|\bm{f}_2-\bm{A}_{12}^T\bm{x}_1^i-\bm{A}_{22}\bm{x}_2^{i,j}\|_{l^2}, \label{eq:innerSenr}
\end{align}
where this time there is no factor $h$ because the condition number of $\bm{A}_{22}$ is bounded. Apart from this modification, the derivation of this estimator essentially follows the same pattern as for \eqref{eq:estimfem}. Once the stopping criterion is reached, after say $j''$ iterations, we set $w_{2,h}^i := w_{2,h}^{i,j''}$ and update $w_h^i = w_{1,h}^i +  w_{2,h}^i$. Again, note that $v_h^i$ is not available in practice, so $\delta_i$ is estimated as \eqref{eq:gsnocheat} instead of \eqref{eq:gscheat}.

The algorithm developed for GFEM \& SGFEM schematically takes the form of Algorithm \ref{alg:gfem}.

\begin{algorithm}[!ht]
 \KwData{$h,\bm{A},\bm{f},k,k'$}
 \KwResult{$w^{i^*}_h,i^*,j^*,j'^*$}
 $\epsilon = h, w_{1,h}^0=w_{2,h}^0=0, e^0=\infty, i=j^*=j'^*=0$\;
 \While{$e^i \geq \epsilon/k$}{
 	$i \leftarrow i+1$\;
  $w_{1,h}^{i,0}=w_{1,h}^{i-1}, e_1^{i,0} = \infty, j=0$\;
 \While{$e_1^{i,j} \geq e^i/k'$}{
 	 $j \leftarrow j+1$\;
     Compute $w_{1,h}^{i,j}$ using initialization $w_{1,h}^{i,j-1}$\;
     Compute error estimator $e_1^{i,j}$ using \eqref{eq:innerSfem}\;
 }
 $j^* \leftarrow j^* + j, w_{1,h}^i = w_{1,h}^{i,j}$\;
 $w_{2,h}^{i,0}=w_{2,h}^{i-1}, e_2^{i,0} = \infty, j=0$\;
 \While{$e_2^{i,j} \geq e^i/k'$}{
 	 $j \leftarrow j+1$\;
     Compute $w_{2,h}^{i,j}$ using initialization $w_{2,h}^{i,j-1}$\;
     Compute error estimator $e_2^{i,j}$ using \eqref{eq:innerSenr}\;
 }
 $j'^* \leftarrow j'^* + j, w_{2,h}^i = w_{2,h}^{i,j}$\;
 $w_h^i = w_{1,h}^i +  w_{2,h}^i$\;
 Compute error estimator $e^i$ using \eqref{eq:gsnocheat}\;
 }
 $i^*=i$.
  \caption{Algorithm for GFEM \& SGFEM.}
 \label{alg:gfem}
\end{algorithm}

Note that, as always, we work on the scaled system \eqref{goodLinSys}. The algorithm stops the outer iterations as soon as the estimated truncation error $e^i$ becomes significantly smaller than the a priori estimated discretization error $\epsilon$. The factor $k$ is there to control the different constants of proportionality appearing in the intermediate calculations (see Appendix \ref{app:stoppingcrit}). During each outer iteration, the algorithm stops the inner iterations as soon as the estimated truncation error $e^{i,j}_\cdot$ becomes significantly lower than the current estimated truncation error $e^i$. Again, the factor $k'$ is there to control the different constants of proportionality appearing in the intermediate calculations (see Appendix \ref{app:stoppingcrit}).

We show in Tables \ref{tab:straight}--\ref{tab:circular} the performances of the iterative solver described above in terms of number of iterations and of computing time when $h$ varies. The FEM has been included in these tables only for the purpose of comparison. For a given $h$, the iterative solver stops when the truncation error tolerance $e^i = \frac{\epsilon}{k}$ with $k=100$ is reached. We used $\epsilon = h^{1/2}$ for the FEM (Algorithm 1) and $\epsilon = h$ for the GFEM/SGFEM (Algorithm 2); these values are based on the a priori discretization error estimates of FEM and GFEM/SGFEM.

For GFEM \& SGFEM, the number of iterations is displayed under the form $i^*\,(j^*,j'^*)$, where $i^*$ is the number of outer iterations, $j^*$ is the cumulated number of CG preconditioned by FMG iterations in $S_{FEM}^h$ and $j'^*$ is the cumulated number of CG iterations in $S_{ENR}^h$. We used a computer with 64 bits architecture, a 3.6GHz processor and 16 GB of RAM. On this computer, we measured that Matlab was able to ``count" up to $150 \times 10^6$ in 1 second. The parameters for these simulations were $a_0 = 1$, $a_1 = 10$, $k = 100$ and $k' = 4$. We also mention that the relaxation scheme used in FMG was the Gauss-Seidel method and the associated stopping criterion was to stop relaxing when the $l^2$-norm of the residual was higher than half of the previous one. Finally, the coarsest level in the FMG scheme is associated with the uniform mesh $\mc T_h$ with $h = 1$, and subsequent levels are given by halving the mesh size $h$, and all initializations are done with the zero vector.

Table \ref{tab:straight} illustrates the case of the straight interface with $\theta_0 = \pi/6$ and $d_0 = 1-1/\sqrt{2}$.

\begin{table}[!ht]
\caption{Results of the CG preconditioned by FMG iterative solver for FEM/GFEM/SGFEM on the straight interface problem with $\theta_0 = \pi/6$.}
\centering
\footnotesize
\begin{tabular}{c|cc|cc|cc}
$1/h$ & \multicolumn{2}{c|}{FEM} & \multicolumn{2}{c|}{GFEM} & \multicolumn{2}{c}{SGFEM} \\ \hline
& \# it. & t (s) & \# it. & t (s) & \# it. & t (s) \\
2 & 2 & 0.0009 & 38 (38,59) & 0.0233  & 8 (8,8) & 0.0047 \\
4 & 2 & 0.0028 & 46 (46,70) & 0.0566  & 12 (12,12) & 0.0174 \\
8 & 2 & 0.0039 & 53 (53,81) & 0.1109  & 12 (12,12) & 0.0301  \\
16 & 2 & 0.0078 & 59 (59,89) & 0.2035 & 13 (13,13) & 0.0449  \\
32 & 3 & 0.0151 & 67 (67,103) & 0.3637 & 13 (13,13) & 0.0710  \\
64 & 3 & 0.0268 & 77 (77,118)  & 0.7558 & 14 (14,14) & 0.1393  \\
128 & 3 & 0.0609 & 83 (113,128) & 2.614 & 14 (14,14) & 0.3374  \\
256 & 4 & 0.2618 & 90 (141,139) & 11.12  & 14 (14,14) & 1.130  \\
512 & 4 & 1.119 &  97 (167,150) & 54.57   & 15 (28,15) & 9.048 \\
1024 & 5 & 5.509  & 103 (176,159) & 227.6  & 15 (28,15) & 35.91 \\
2048 & 6 & 30.60 & 112 (229,173) & 1302 & 15 (28,15) & 158.6
\end{tabular}
\normalsize
\label{tab:straight}
\end{table}

Table \ref{tab:patho} illustrates the case of the straight interface with $\theta_0 = \pi/4$ (pathological case where the interface is parallel to some of the mesh edges) and $d_0$ such that the relative distance between the mesh and the interface was only $10^{-3}$. 

\begin{table}[!ht]
\caption{Results of the CG preconditioned by FMG iterative solver for FEM/GFEM/SGFEM on the straight interface problem with $\theta_0 = \pi/4$ and relative distance to the mesh $10^{-3}$.}
\centering
\footnotesize
\begin{tabular}{c|cc|cc|cc}
$1/h$ & \multicolumn{2}{c|}{FEM} & \multicolumn{2}{c|}{GFEM} & \multicolumn{2}{c}{SGFEM} \\ \hline
& \# it. & t (s) & \# it. & t (s) & \# it. & t (s) \\
2 & 3 & 0.0017 & 44 (54,69) & 0.0362 & 13 (13,13) & 0.0088 \\
4 & 3 & 0.0039 & 33 (34,49) & 0.0432 & 17 (18,17) & 0.0210 \\
8 & 3 & 0.0056 & 5 (6,5) & 0.0120 & 5 (6,5) & 0.0119 \\
16 & 3 &  0.0092 & 5 (6,5) & 0.0193 & 5 (6,5) & 0.0233 \\
32 & 5 & 0.0241 & 5 (6,5) & 0.0305 & 5 (6,5) & 0.0305 \\
64 & 5 & 0.0423 & 16 (21,23) & 0.1976 & 7 (10,7) & 0.0915 \\
128 & 5 & 0.0955 & 128 (166,199) & 3.772 & 7 (10,7) & 0.2131 \\
256 & 6 & 0.3767 & 119 (170,185) & 13.13 & 9 (12,9) & 0.9202 \\
512 & 6 & 1.565 & 84 (133,130) & 42.04 & 11 (14,11) & 4.334 \\
1024 & 7 & 7.609 & 53 (90.69) & 126.1 & 18 (22,18) & 30.62 \\
2048 & 9 & 42.91 & 140 (571,218) & 3128 & 7 (15,7) & 74.25
\end{tabular}
\normalsize
\label{tab:patho}
\end{table}

Table \ref{tab:circular} illustrates the case of the circular interface with $r_c=1/\sqrt{10}$ and $(x_c,y_c) = (1/\sqrt{5},1/\sqrt{3})$.

\begin{table}[!ht]
\caption{Results of the CG preconditioned by FMG iterative solver for FEM/GFEM/SGFEM on the circular interface problem.}
\centering
\footnotesize
\begin{tabular}{c|cc|cc|cc}
$1/h$ & \multicolumn{2}{c|}{FEM} & \multicolumn{2}{c|}{GFEM} & \multicolumn{2}{c}{SGFEM} \\ \hline
& \# it. & t (s) & \# it. & t (s) & \# it. & t (s) \\
2 & 2 & 0.0009 & 6 (6,6) & 0.0036  & 6 (6,6) & 0.0035 \\
4 & 2 & 0.0021 & 22 (22,22) & 0.0255  & 8 (8,8) & 0.0093 \\
8 & 3 & 0.0057 & 23 (23,23) & 0.0478  & 10 (10,10) & 0.0200  \\
16 & 3 & 0.0092 & 49 (49,64) & 0.1681 & 10 (10,10) & 0.0326  \\
32 & 4 & 0.0192 & 62 (63,94) & 0.3506 & 13 (13,13) & 0.0667  \\
64 & 4 & 0.0344 &  83 (115,131) & 1.109 & 16 (17,16) & 0.1610  \\
128 & 5 & 0.0970 & 92 (131,147) & 2.940 & 15 (18,15) & 0.3932  \\
256 & 6 & 0.3669 & 100 (165,166) & 11.25  & 16 (30,16) & 1.995  \\
512 & 6 & 1.445 &  120 (271,201) & 74.48   & 17 (32,17) & 9.169 \\
1024 & 7 & 7.217  & 142 (338,244) & 398.2  & 20 (40,20) & 47.01 \\
2048 & 8 & 37.41 & 156 (387,270) & 1976  & 22 (51,22) & 259.4
\end{tabular}
\normalsize
\label{tab:circular}
\end{table}

The \# it. and t(s) for FEM for a given $h$ in Tables \ref{tab:straight}--\ref{tab:circular} are much less than for GFEM/SGFEM. This is expected as the error in FEM for a given $h$ ($O(h^{1/2})$) is much greater than that of GFEM/SGFEM ($O(h)$). Moreover, we observe that the computational time $t$ scales a little over quadratically with the mesh size and we have roughly $t = O(h^{-2.3})$ (except for the pathological case of the straight interface problem with $\theta_0 = \pi/4$ and ``small" relative distance to the mesh). This rate is slightly over the optimal rate of $O(h^{-2})$ because although we use efficient solvers, our stopping criteria are somewhat pessimistic since they rely on the norm of the residual. Concerning the pathological case of the interface parallel to the mesh edges, Table \ref{tab:patho} reveals that GFEM is not stable, while SGFEM is. Remember that in Figures \ref{fig:straightcondd0} and \ref{fig:straightangled0} we saw that the condition number of (M-)GFEM blew up as the interface was getting closer and closer to the mesh edges, while the angle was going to 0. Conversely, SGFEM was stable in this situation.

We also note that for the three situations shown in Tables \ref{tab:straight}-\ref{tab:circular}, many more outer iterations are needed for GFEM than for SGFEM. This is a direct consequence of the fact that the angle for SGFEM is larger than that for the GFEM, as indicated in Figures \ref{fig:straightangleh} and \ref{fig:circleangleh}. Indeed, if we were to solve the block Gauss-Seidel system \eqref{eq:Sfem}--\eqref{eq:Senr} exactly (i.e., we only had outer iterations) and if there had been only one angle between the spaces (i.e., the smallest) then, for the same decrease in the truncation error, for each outer iteration in SGFEM we would have needed $n = \frac{\log \left( \cos \vartheta_2 \right)}{\log \left( \cos \vartheta_1 \right)}$ outer iterations in GFEM, where $\vartheta_1$ denotes the angle for GFEM and $\vartheta_2$ denotes the angle for SGFEM. Indeed, each iteration will result in a reduction of the truncation error by a factor $q_1 = \cos^2 \vartheta_1$ for GFEM and $q_2 = \cos^2 \vartheta_2$ for SGFEM. As a result, $n$ iterations for GFEM will decrease the error by $q_1^n$. Solving $q_1^n = q_2$ for $n$ yields the above ratio in terms of angles $\vartheta_1,\vartheta_2$. We have in fact solved the block Gauss-Seidel system \eqref{eq:Sfem}--\eqref{eq:Senr} by using direct solvers for the blocks and the stopping criterion \eqref{eq:gscheat}. We did observe the role of the angle, i.e., SGFEM yielded a speed-up of roughly 6-9 times compared to the GFEM. However our iterative solver with two estimators performed much faster (by a factor of about 7 for the values of $h$ considered in our experiments) on both GFEM and SGFEM than solving the block Gauss-Seidel system \eqref{eq:Sfem}--\eqref{eq:Senr} directly. We do not show the results as in this paper we are only concerned with our iterative scheme. Note that our error estimators estimate the truncation error as well as the discretization error and then ``balance" these two errors; these estimators are very fast to compute.

However, since we do not solve \eqref{eq:Sfem}--\eqref{eq:Senr} exactly and use estimates to decide when to stop the inner iterations, this ratio $n$ mentioned in the last paragraph is slightly perturbed. Recall also that the angle we discuss in this paper is in fact the smallest between the spaces. There are other, larger, angles that could affect the estimation of the truncation error since we use Richardson extrapolation to estimate $\delta_i$ using \eqref{eq:gsnocheat}. Since the angle for SGFEM is larger than for GFEM, this ``pollution" by larger angles affects GFEM more than SGFEM. As a result, the risk of under-resolving system \eqref{goodLinSys} is higher for GFEM than SGFEM.

We also considered other solvers not shown here (e.g., V-cycles as a solver, FMG followed by V-cycles as a solver, CG preconditioned by some V-cycles), however, the presented solver -- CG preconditioned by FMG -- was found to be the most robust and computationally efficient. The reason is that CG preconditioned by multigrid is more robust than multigrid alone, and that FMG is more efficient than V-cycles.

Let us verify if the iterative solutions have converged to discretization accuracy. We now denote by $v_h \in S^h$ the iterated solution of \eqref{goodLinSys} yielded by the chosen iterative solver and by $\hat{\epsilon}^h = \| u - v_h\|_{\mc E}$ the total error due to both discretization and truncation. Thanks to Galerkin orthogonality between $u - u_h$ and $u_h - v_h$ in the $B(\cdot,\cdot)$ inner product, it holds, for the straight interface problem
\begin{align*}
\left( \hat{\epsilon}^h \right)^2 &= \| u  - v_h \|^2_{\mc E},\\
 &= \| u  - u_h \|^2_{\mc E} + \| u_h  - v_h \|^2_{\mc E},\\
 &= \left( \epsilon^h \right)^2 + \delta^2,
\end{align*}
so that the total error $\hat{\epsilon}^h$ can be orthogonally decomposed into discretization error $\epsilon^h$ and truncation error $\delta$. Similarly, for the circular interface problem we set $\left( \hat{\epsilon}^h \right)^2 = \left( \epsilon^h \right)^2 + \delta^2$, so that
\begin{align*}
\left( \hat{\epsilon}^h \right)^2 &= \left|\| u\|^2_{\mc E} - \| u_h\|^2_{\mc {\widetilde{E}}} \right| + \| u_h  - v_h \|^2_{\mc {\widetilde{E}}}.
\end{align*}
In practice, we also solve \eqref{goodLinSys} using a direct solver in order to have $u_h$. We can thus compute $\epsilon^h$ as well as $\hat{\epsilon}^h$.

\begin{figure}[!ht]
\centerline{\includegraphics[scale = 0.6]{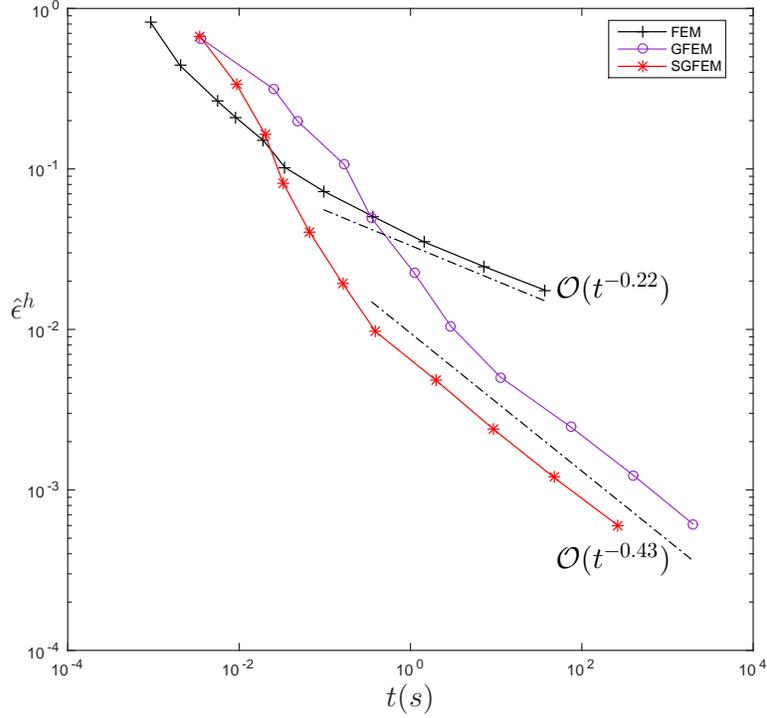}}
\caption{Evolution of the (relative) error $\hat{\epsilon}^h$ as the computational time varies on the circular interface problem.}
\label{fig:ehvsetotcircular}
\end{figure}

Figure \ref{fig:ehvsetotcircular} displays the evolution of the (relative) total error $\hat{\epsilon}^h$ as the computational time $t$ varies, on the circular interface problem (same parameters as for Table \ref{tab:circular}). We observe that the error $\hat{\epsilon}^h$ scales with the computational time $t$ for FEM as $\hat{\epsilon}^h =0.04\times t^{-0.22}$. For GFEM, we have roughly $\hat{\epsilon}^h = 0.015\times t^{-0.43}$. For SGFEM, we have $\hat{\epsilon}^h = 0.006\times t^{-0.43}$. The optimal rates would be $O(t^{-1/4})$ for FEM and $O(t^{-1/2})$ for GFEM \& SGFEM. Moreover, the SGFEM is roughly eight or nine times faster than the GFEM for the same accuracy. As discussed earlier, this is a direct consequence of the angle property. We also mention that for the straight interface problem (with the same parameters as for Table \ref{tab:straight}) for FEM we have $\hat{\epsilon}^h =0.0275\times t^{-0.22}$. For GFEM, we have roughly $\hat{\epsilon}^h = 0.0075\times t^{-0.43}$. For SGFEM, we have $\hat{\epsilon}^h = 0.0035\times t^{-0.43}$. This in turn tells us that the SGFEM is roughly six times faster than the GFEM for the same accuracy on this problem.

\begin{figure}[!ht]
\centerline{\includegraphics[scale = 0.65]{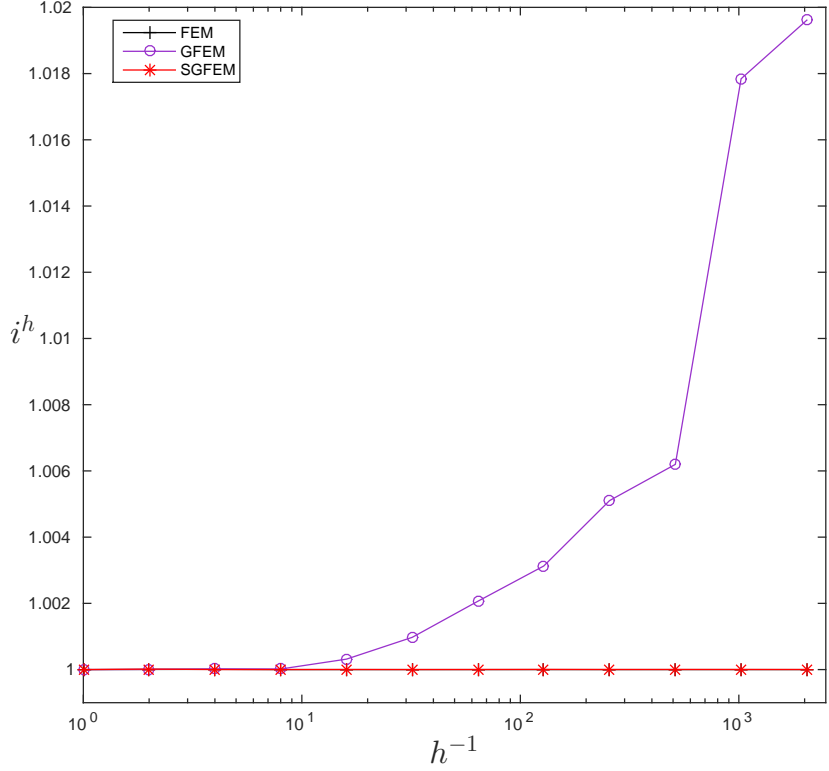}}
\caption{Evolution of the ratio iterative solution error over discretization error against $h$.}
\label{fig:ehvsetot}
\end{figure}

To verify if our iterative solutions have indeed converged to discretization accuracy, we show in Figure \ref{fig:ehvsetot} the evolution of the ratio iterative solution error over discretization error $i^h := \displaystyle \frac{\hat{\epsilon}^h}{\epsilon^h} = \sqrt{1+\left(\frac{\delta}{\epsilon^h}\right)^2}$ (on the circular interface problem with the same parameters as before). We observe that the ratio $i^h$ is very close to unity for all considered methods. There is no appreciable difference for FEM or SGFEM, while for GFEM the ratio $i^h$ is less than 2\% over unity. Our solutions have thus converged to discretization accuracy. As discussed before, GFEM performs worse than SGFEM because of the angle property and the pollution by larger angles in the estimation of the truncation error.

\begin{figure}[!ht]
\centerline{\includegraphics[scale = 0.65]{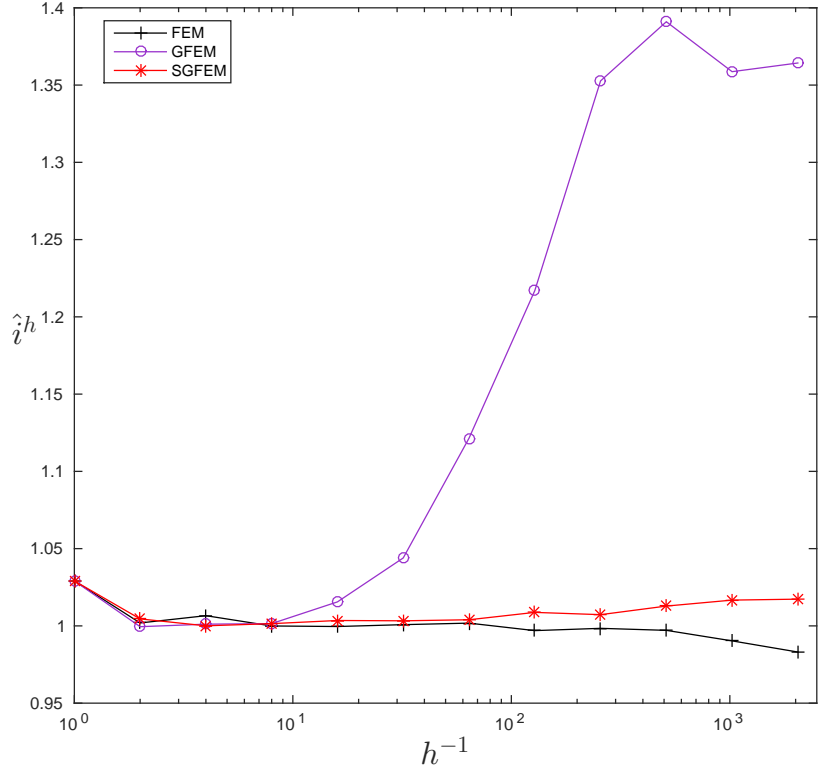}}
\caption{Efficiency of the error estimator $\hat{i^h}$ as $h$ varies.}
\label{fig:ivsh}
\end{figure}

A last remark has to be made about the solutions yielded by the iterative solvers in the case where the exact solution $u \in \mc E$ and the discrete solution $u_h \in S^h$ are unknown. Once the iterative solutions have been computed they can be used to design an error estimator. Indeed, neither the solver nor the stopping criteria rely on the knowledge of $u$ or $u_h$. Using Richardson extrapolation this time on the norm of $v_h$ allows us to extrapolate when $h \rightarrow 0$ to obtain some $\eta$ approximating $\|u\|_\mc{E}$. Indeed
\begin{align*}
\|u\|_\mc{E} &= \|v_h\|_\mc{E} + \varepsilon_h,\\
\|u\|_\mc{E} &= \|v_{2h}\|_\mc{E} + \varepsilon_{2h}.
\end{align*}
Next, assuming $\varepsilon_h = Ch^p$ where $p$ is known by the underlying properties of the PDE, the choice of the partition of unity and the choice of the enrichment, we have
\begin{align*}
C = \frac{\|v_h\|_\mc{E}  - \|v_{2h}\|_\mc{E} }{(2^p-1)h^p}.
\end{align*}
As a result, we can estimate $\|u\|_\mc{E}$ as
\begin{align*}
\eta = \|v_h\|_\mc{E} + Ch^p.
\end{align*}
Finally, by assuming Galerkin orthogonality, we have an error estimator using $\bar{\epsilon}^h := \left| \eta^2 - \|v_h\|_\mc{E}^2 \right|^{1/2}$. The graph of $\bar{\epsilon}^h$ as the computational time $t$ varies is similar to Figure \ref{fig:ehvsetotcircular} and we do not include it here. The efficiency of this estimator can then be assessed by computing $\hat{i}^h  := \displaystyle\frac{\bar{\epsilon}^h}{\hat{\epsilon}^h}$. The results are presented in Figure \ref{fig:ivsh} (on the circular interface problem with the same parameters as before). We observe that the efficiency of the error estimator $\hat{i}^h$ stays close to unity for all considered methods: $\hat{i}^h$ is within the 3\% range for FEM and SGFEM, 40\% for GFEM. Once again, GFEM performs worse than SGFEM because of the angle property and the pollution by larger angles in the estimation of the truncation error.

We summarize the results mentioned above about the differences between GFEM and SGFEM in the following table, where for a given relative error tolerance $\tau = 5\%, 1\%, 0.1\%$, we show the computed relative total error $\hat{\epsilon}^h$, the computed relative discretization error $\epsilon^h$, the computed relative truncation error $\delta$, the computed relative extrapolated error $\bar{\epsilon}^h$, time $t(s)$ and the efficiency indexes $i^h$ and $\hat{i}^h$, on the circular interface problem with the same parameters as before.

\begin{center}
\begin{tabular}{c|ccccccc}
 & \multicolumn{6}{c}{GFEM} \\ \hline
 $\tau$ & $\hat{\epsilon}^h$ & $\epsilon^h$ & $\delta$ & $\bar{\epsilon}^h$ & $t(s)$ & $i^h$ & $\hat{i}^h$ \\ \hline
5\% & 4.95\% & 4.94\% & 0.22\% & 5.16\% & 0.3506 & 1.001 & 1.044 \\ 
1\% & 1.05\% & 1.05\% & 0.083\% & 1.28\% & 2.940 & 1.003 & 1.217 \\ 
0.1\% & 0.123\% &    0.121\% &    0.023\% &0.168\% & 398.2 & 1.018 & 1.359
\end{tabular}
\end{center}

\begin{center}
\begin{tabular}{c|ccccccc}
 & \multicolumn{6}{c}{SGFEM} \\ \hline
 $\tau$ & $\hat{\epsilon}^h$ & $\epsilon^h$ & $\delta$ & $\bar{\epsilon}^h$ & $t(s)$ & $i^h$ & $\hat{i}^h$ \\ \hline
5\% & 4.03\% &  4.03\% &  0.0069\% & 4.05\%  & 0.0667 & 1.000 & 1.003 \\ 
1\% & 0.967\% &  0.967\% &  0.0032\% & 0.976\% & 0.3932 & 1.000 & 1.009 \\ 
0.1\% & 0.120\% &  0.120\% &  0.00032\% & 0.122\%  & 47.01 & 1.000 & 1.017
\end{tabular}
\end{center}

\section*{Conclusion}

In this paper, we have considered several GFEMs on an interface problem. The proposed GFEMs differed with regards to the choice of enrichment function and space. We have illustrated how these different choices could allow recovering the optimal order of convergence of the underlying partition of unity. However, as we have shown, this so-called approximation property is not the only important feature of a well-designed GFEM, and it is equally important that the GFEM be well-conditioned in order to be able to solve the resulting linear system efficiently. To this end, we have highlighted the importance of the so-called angle between the approximation spaces. We emphasized how this angle was related to the conditioning of GFEM and to the stability of the method with respect to varying criteria, such as the position of the interface with respect to the mesh. Finally, we have showed in the last section that this angle property could be exploited within a block Gauss-Seidel iterative scheme between the approximation spaces, resulting in large savings in computational times.

\appendix
\section{Angle between subspaces}
\label{app:angle}
In this first annex, we derive the formula yielding the angle $\vartheta$ between subspaces $S_1$ and $S_2$, of dimensions $m$ and $n$ respectively, in the sense of the inner product $B(\cdot,\cdot)$.\\
Let $S = S_1 +  S_2 = \left\lbrace s = (s_1,s_2) : s_1\in S_1, s_2\in S_2 \right\rbrace$ and $P_1$ be the orthogonal projection operator on $S_1$. Then the angle $\vartheta$ is defined by
\begin{equation*}
\cos \vartheta = \max\limits_{s_2 \in S_2} \displaystyle \frac{\|P_1(s_2)\|_\mathcal{E}}{\|s_2\|_\mathcal{E}}.
\end{equation*}
By definition of the projection $P_1$, we have
\begin{equation}
B(P_1(s_2),s_1) = B(s_2,s_1), \quad\forall s_1\in S_1.
\label{proj}
\end{equation}
Identifying $s_1 \approx X_1 \in \mathbb{R}^m$, $s_2 \approx X_2 \in \mathbb{R}^n$ and $P_1(s_2) \approx W \in \mathbb{R}^m$, we can then use the matrix $\bm{A}$ defined in \eqref{goodLinSys}, so problem \eqref{proj} reads
\begin{equation*}
(W^T,0)\bm{A}(X_1,0)^T = (0,X_2^T)\bm{A}(X_1,0)^T,\quad \forall X_1\in \mathbb{R}^m.
\end{equation*}
With the sub-matrices defined in \eqref{goodLinSys}, we have
\begin{equation*}
 W^T\bm{A}_{11}X_1 = X_2^T \bm{A}_{21} X_1,\quad \forall X_1\in \mathbb{R}^m
\end{equation*}
which leads to ($\bm{A}$ being symmetric)
\begin{equation*}
 W = \bm{A}_{11}^{-1} \bm{A}_{12} X_2.
\end{equation*}
Then, one can simply derive
\begin{align*}
\|P_1(s_2)\|^2_\mathcal{E} &= W^T \bm{A}_{11} W,\\ 
&= X_2^T\bm{A}_{21}\bm{A}_{11}^{-1}\bm{A}_{12}X_2,
\end{align*}
and we obviously have
\begin{equation*}
\|s_2\|^2_\mathcal{E} = X_2^T\bm{A}_{22}X_2.
\end{equation*}
Hence
\begin{equation}
\cos^2 \vartheta = \max\limits_{X_2 \in \mathbb{R}^n} \displaystyle \frac{X_2^T\bm{A}_{21}\bm{A}_{11}^{-1}\bm{A}_{12}X_2}{ X_2^T\bm{A}_{22}X_2},
\label{ans}
\end{equation}
which leads to the generalized eigenvalue problem
\begin{equation*}
\bm{A}_{21}\bm{A}_{11}^{-1}\bm{A}_{12} = \lambda \bm{A}_{22}.
\end{equation*}
The largest eigenvalue of this problem is equal to $\max\limits_{X_2 \in \mathbb{R}^n}  \frac{X_2^T\bm{A}_{21}\bm{A}_{11}^{-1}\bm{A}_{12}X_2}{ X_2^T\bm{A}_{22}X_2}$, so one can then find $\vartheta$ using \eqref{ans}.

\section{Error induced by the perturbation}
\label{app:error}
In this second annex, we show 
\begin{align*}
\left| \| u -u_h \|^2_{\mc E} -  \left|\|u\|^2_\mc{E}-\|u_h\|^2_\mc{\widetilde{E}}\right| \right| &\leq O(h^k),
\end{align*}
where $k = 3/2$ for FEM and $k = 2$ for GFEM \& SGFEM, $u \in \mc E$ is the (exact) solution of \eqref{InterfaceProb1d}, $u_h \in S^h$ is the (discrete perturbed) solution of \eqref{Interfaceperturbfinitedim}, $\|\cdot\|_{\mc E}:= B(\cdot,\cdot)^{1/2}$ is the usual energy norm and $\|\cdot\|_{\mc {\widetilde{E}}}:= \widetilde{B}(\cdot,\cdot)^{1/2}$ is the perturbed energy norm. We also introduce $\widetilde{u} \in \mc E$ the (exact perturbed) solution of \eqref{Interfaceperturb}.

The proof is divided into two parts: in the first part, we show that
$\left| \|u\|^2_{\mc E} - \|\widetilde{u}\|^2_{\mc {\widetilde{E}}}\right| = O(h^2)$. In the second part, we use this result to show that $\left| \| u -u_h \|^2_{\mc E} -  \left|\|u\|^2_\mc{E}-\|u_h\|^2_\mc{\widetilde{E}}\right| \right| \leq O(h^k)$, where $k = 3/2$ for FEM and $k = 2$ for GFEM \& SGFEM.

We have two weak formulations: one for the original problem \eqref{InterfaceProb1d} and one for the perturbed problem \eqref{Interfaceperturb}
\begin{align}
B(u,v) =& \int_\Omega a \nabla u \cdot \nabla v \, d\bm{x} = \int_{\partial\Omega} g_N v\, ds, \quad \forall v \in \mc{E},\label{weakB}\\
\widetilde{B}(\widetilde{u},v) =& \int_\Omega \widetilde{a} \nabla \widetilde{u}\cdot \nabla v \, d\bm{x} = \int_{\partial\Omega} g_N v\, ds, \quad \forall v \in \mc{E}. \label{weakBtilde}
\end{align}
Let us first show that the two norms induced by $B(\cdot,\cdot)$ and $\widetilde{B}(\cdot,\cdot)$ are equivalent on $\mc{E}$. For all $v \in \mc{E}$, we have
\begin{align}
\beta_0 \int_\Omega \nabla v \cdot \nabla v \, d\bm{x} &\leq B(v,v) \leq \beta_1 \int_\Omega \nabla v \cdot \nabla v \, d\bm{x},\nonumber\\
\beta_0 | v |_{H^1}^2 &\leq B(v,v) \leq \beta_1 | v |_{H^1}^2, \label{normbounds}
\end{align}
where $\beta_0,\beta_1$ were defined in Section \ref{sec2} and represent bounds on the coefficient $a$. As a result, the norm induced by $B(\cdot,\cdot)$ and the $H^1$ semi-norm are equivalent on $\mc{E}$. The same reasoning follows for $\widetilde{B}(\cdot,\cdot)$, which is also equivalent to $| \cdot |_{H^1}$ on $\mc{E}$. As a result, $B(\cdot,\cdot)$ and $\widetilde{B}(\cdot,\cdot)$ are equivalent on $\mc{E}$. We further emphasis that the bounds appearing in \eqref{normbounds} are independent of $h$. It follows that $B(v,v) = O(h^p)$ holds if and only if $\widetilde{B}(v,v) = O(h^p)$ holds.

Now, considering \eqref{weakB} and \eqref{weakBtilde}, we have
\begin{align}
B(u,v) &= \widetilde{B}(\widetilde{u},v),\quad \forall v \in \mc{E}.\label{BBtilde0}
\end{align}
We immediately obtain
\begin{align}
B(u,u) = \widetilde{B}(\widetilde{u},u), \mbox{ and } B(u,\widetilde{u}) &=  \widetilde{B}(\widetilde{u},\widetilde{u}).\label{BBtilde1}
\end{align}

We also have
\begin{align}
B(u-\widetilde{u},v) &= B(u,v)-B(\widetilde{u},v), \mbox{ and, using } \eqref{BBtilde0},\nonumber\\
&= \widetilde{B}(\widetilde{u},v) - B(\widetilde{u},v),\nonumber\\
&= \int_\Omega (\widetilde{a}-a) \nabla \widetilde{u}\cdot \nabla v \, d\bm{x},\nonumber\\
&= (a_0 - a_1) \int_\omega \nabla \widetilde{u}\cdot \nabla v \, d\bm{x},\quad \forall v \in \mc{E}. \label{precauchy1}
\end{align}
Using Cauchy-Schwarz inequality, it yields
\begin{align}
\left|B(u-\widetilde{u},v)\right| &\leq  \left|a_0-a_1\right|\left(\int_\omega \nabla \widetilde{u}\cdot \nabla \widetilde{u} \, d\bm{x} \right)^{1/2} \left(\int_\omega \nabla v \cdot \nabla v \, d\bm{x} \right)^{1/2},\nonumber\\
&\leq  \left|a_0-a_1\right| \| \nabla \widetilde{u} \|_{L^2(\omega)} \| \nabla v \|_{L^2(\omega)}. \label{cauchy1}
\end{align}
From this we can obtain
\begin{align}
\left|B(u,u) - \widetilde{B}(\widetilde{u},\widetilde{u})\right| &= \left|B(u,u) - B(u,\widetilde{u})\right|, \mbox{ using } \eqref{BBtilde1},\nonumber\\
&= \left|B(u,u-\widetilde{u})\right|,\nonumber\\
&\leq \left|a_0-a_1\right| \| \nabla \widetilde{u} \|_{L^2(\omega)} \| \nabla u \|_{L^2(\omega)}, \mbox{ using } \eqref{cauchy1} \mbox{ with } v=u,\nonumber\\
&\leq  \left|a_0-a_1\right| \mu(\omega)^{1/2} \| \nabla \widetilde{u} \|_{L^2(\omega)} \| \nabla u \|_{L^\infty(\omega)},\label{keyeq}
\end{align}
where we have used $\| \nabla u \|_{L^2(\omega)} \leq  \mu(\omega)^{1/2} \| \nabla u \|_{L^\infty(\omega)}$ and the term on the right hand side is finite because our manufactured solution $u$ does not have any singularity in $\ov{\Omega}$. Now, using $\mu(\omega) = O(h^2)$, we get
\begin{align}
\left|B(u,u) - \widetilde{B}(\widetilde{u},\widetilde{u})\right| &= O(h). \label{orderfem}
\end{align}

Similarly to \eqref{precauchy1}
\begin{align}
\widetilde{B}(u-\widetilde{u},v) &= \widetilde{B}(u,v) - \widetilde{B}(\widetilde{u},v), \mbox{ and, using } \eqref{BBtilde0},\nonumber\\
&= \widetilde{B}(u,v) - B(u,v),\nonumber\\
&= (a_0-a_1)\int_\omega \nabla u\cdot \nabla v \, d\bm{x},\quad \forall v \in \mc{E}. \label{precauchy2}
\end{align}
Let us now look at
\begin{align*}
\widetilde{B}(u-\widetilde{u},u-\widetilde{u}) &= \widetilde{B}(u-\widetilde{u},u) - \widetilde{B}(u-\widetilde{u},\widetilde{u}),\\
&= (a_0-a_1)\| \nabla u \|^2_{L^2(\omega)} - \widetilde{B}(u-\widetilde{u},\widetilde{u}), \mbox{ using } \eqref{precauchy2} \mbox{ with } v=u,\\
&= (a_0-a_1)\| \nabla u \|^2_{L^2(\omega)} +\widetilde{B}(\widetilde{u},\widetilde{u}) -\widetilde{B}(u,\widetilde{u}),\\
&= (a_0-a_1)\| \nabla u \|^2_{L^2(\omega)} +\widetilde{B}(\widetilde{u},\widetilde{u}) -B(u,u), \mbox{ using } \eqref{BBtilde1}.
\end{align*}
Thus
\begin{align}
\widetilde{B}(u-\widetilde{u},u-\widetilde{u}) &\leq \left|a_0-a_1\right|\| \nabla u \|^2_{L^2(\omega)} +\left|B(u,u)-\widetilde{B}(\widetilde{u},\widetilde{u})\right|, \mbox{ using triangular inequality},\nonumber\\
&\leq \left|B(u,u)-\widetilde{B}(\widetilde{u},\widetilde{u})\right| + O(h^2), \label{pivot}
\end{align}
where again we have used $\| \nabla u \|_{L^2(\omega)} \leq  \mu(\omega)^{1/2} \| \nabla u \|_{L^\infty(\omega)}$ and $\mu(\omega) = O(h^2)$.

We now have most of the ingredients to prove that $\left|B(u,u)-\widetilde{B}(\widetilde{u},\widetilde{u})\right| = O(h^2)$. We will proceed by induction. Let us consider the two statements
\begin{align*}
\left\lbrace\begin{array}{ll}
\left|B(u,u)-\widetilde{B}(\widetilde{u},\widetilde{u})\right| = O(h).\\
\left|B(u,u)-\widetilde{B}(\widetilde{u},\widetilde{u})\right| = O(h^p) \Rightarrow \left|B(u,u)-\widetilde{B}(\widetilde{u},\widetilde{u})\right| = O(h^{1+p/2}).
\end{array} \right.
\end{align*}

The first statement has already been proven in \eqref{orderfem}. Let us prove the second. Assume that $\left|B(u,u)-\widetilde{B}(\widetilde{u},\widetilde{u})\right| = O(h^p)$ holds for some $1 \leq p \leq 2$. Then, using the induction assumption in \eqref{pivot}, we have
\begin{align*}
\widetilde{B}(u-\widetilde{u},u-\widetilde{u}) &\leq O(h^p) + O(h^2),\\
&\leq O(h^p).
\end{align*}

Then, by the equivalence of $B(\cdot,\cdot)$ and $\widetilde{B}(\cdot,\cdot)$ on $\mc{E}$, we also have
\begin{align}
B(u-\widetilde{u},u-\widetilde{u}) &\leq O(h^p). \label{eq:equiv}
\end{align}
Now, consider the following
\begin{align}
\left|B(u-\widetilde{u},u)\right| &= \left|B(u,u) - B(\widetilde{u},u)\right|,\nonumber\\
&= \left|B(u,u) - \widetilde{B}(\widetilde{u},\widetilde{u})\right|, \mbox{ using } \eqref{BBtilde1},\nonumber\\
&\leq O(h^p), \mbox{ using the induction assumption again}. \label{eq:equiv2}
\end{align}
Now, using \eqref{precauchy1} with $v = \widetilde{u}$
\begin{align*}
\left|a_0-a_1\right| \| \nabla \widetilde{u} \|^2_{L^2(\omega)} & = \left|B(u-\widetilde{u},\widetilde{u})\right|,\\
&= \left|B(u-\widetilde{u},u) - B(u-\widetilde{u},u-\widetilde{u})\right|,\\
&\leq \left|B(u-\widetilde{u},u)\right| + B(u-\widetilde{u},u-\widetilde{u}), \mbox{ by the triangular inequality},\\
&\leq O(h^p), \mbox{ using \eqref{eq:equiv} and \eqref{eq:equiv2}}.
\end{align*}
Which yields
\begin{align}
 \| \nabla \widetilde{u} \|_{L^2(\omega)} &\leq O(h^{p/2}). \label{eq:0}
\end{align}
Recall \eqref{keyeq}
\begin{align*}
\left|B(u,u) - \widetilde{B}(\widetilde{u},\widetilde{u})\right|
&\leq  \left|a_0-a_1\right| \mu(\omega)^{1/2} \| \nabla \widetilde{u} \|_{L^2(\omega)} \| \nabla u \|_{L^\infty(\omega)},\\
&\leq O(h^{1+p/2}), \mbox{ using \eqref{eq:0} and $\mu(\omega) = O(h^2)$,}
\end{align*}
which is the desired result. Applying it inductively starting at $p=1$ yields
$\left| B(u,u)-\widetilde{B}(\widetilde{u},\widetilde{u})\right| = O(h^2)$. Equivalently, we have $\left| \|u\|^2_{\mc E} - \|\widetilde{u}\|^2_{\mc {\widetilde{E}}}\right| = O(h^2)$. In particular, all the relations written in the induction proof hold for $p=2$.

Now, let us prove the second result, which is
\begin{align*}
\left| \| u -u_h \|^2_{\mc E} -  \left|\|u\|^2_\mc{E}-\|u_h\|^2_\mc{\widetilde{E}}\right| \right| \leq O(h^k),
\end{align*}
where $k = 3/2$ for FEM and $k = 2$ for GFEM \& SGFEM and $u_h\in S^h$ is the solution of the following variational problem
\begin{align*}
\widetilde{B}(u_h,v) = F(v),\quad \forall v \in S^h.
\end{align*}
We split the result into these two inequalities for the sake of clarity: we need to show that
\begin{align*}
 \| u -u_h \|^2_{\mc E} &\leq   \left|\|u\|^2_\mc{E}-\|u_h\|^2_\mc{\widetilde{E}}\right| + O(h^k),\mbox{ and,}\\
\left|\|u\|^2_\mc{E}-\|u_h\|^2_\mc{\widetilde{E}}\right|  &\leq \| u -u_h \|^2_{\mc E}   + O(h^k).
\end{align*}

Using the triangular inequality, it holds
\begin{align}
\| u -u_h \|_{\mc E} &\leq \| u - \widetilde{u} \|_{\mc E} + \| \widetilde{u} -u_h \|_{\mc E},\nonumber\\
&\leq  \| \widetilde{u} -u_h \|_{\mc E} + O(h), \mbox{ using \eqref{eq:equiv} with } p=2. \label{eq:1}
\end{align}

Next, let us consider the difference
\begin{align*}
\left|\| \widetilde{u} - u_h \|^2_{\mc E} - \| \widetilde{u} - u_h \|^2_{\mc{\widetilde{E}}}\right| &= \left|\int_\Omega (a - \widetilde{a}) \nabla (\widetilde{u} - u_h) \cdot \nabla (\widetilde{u} - u_h) \,d \bm x\right|,\\
&= \left|a_1 - a_0\right| \| \nabla (\widetilde{u} - u_h) \|^2_{L^2(\omega)}.
\end{align*}
Using \eqref{eq:0} with $p=2$ yields $\| \nabla \widetilde{u} \|^2_{L^2(\omega)} \leq O(h^2)$. And $\| \nabla u_h \|^2_{L^2(\omega)} \leq \mu(\omega)\| \nabla u_h \|^2_{L^\infty(\omega)} \leq O(h^2)$ as well since $u_h \in S^h$ does not have any singularity in $\ov \Omega$.  As a result

\begin{align}
\left|\| \widetilde{u} - u_h \|^2_{\mc E} - \| \widetilde{u} - u_h \|^2_{\mc{\widetilde{E}}}\right| \leq O(h^2). \label{eq:3}
\end{align}

Then, thanks to Galerkin orthogonality between $\widetilde{u}-u_h$ and $u_h$ in the $\widetilde{B}(\cdot,\cdot)$ inner-product, we have
\begin{align*}
\| \widetilde{u} -u_h \|_{\mc {\widetilde{E}}}^2 &= \| \widetilde{u} \|_{\mc {\widetilde{E}}}^2 - \| u_h \|_{\mc {\widetilde{E}}}^2.
\end{align*}
As a result, it holds
\begin{align}
\| \widetilde{u} - u_h \|^2_{\mc E} \leq \| \widetilde{u} \|_{\mc {\widetilde{E}}}^2 - \| u_h \|_{\mc {\widetilde{E}}}^2 + O(h^2). \label{eq:4}
\end{align}

By a priori error estimation, we have
\begin{align*}
\| \widetilde{u} - u_h \|_{\mc {\widetilde{E}}}  = O(h^p),
\end{align*}
where $p = 1/2$ for FEM and $p = 1$ for GFEM \& SGFEM. By equivalence of $B(\cdot,\cdot)$ and $\widetilde{B}(\cdot,\cdot)$ on $\mc{E}$, it holds
\begin{align}
\| \widetilde{u} - u_h \|_{\mc {E}}  = O(h^p). \label{eq:2}
\end{align}

Hence, starting with \eqref{eq:1}
\begin{align*}
\| u -u_h \|^2_{\mc E} &\leq \left(\| \widetilde{u} -u_h \|_{\mc E} + O(h)\right)^2,\\
&\leq \| \widetilde{u} -u_h \|^2_{\mc E} + \| \widetilde{u} -u_h \|_{\mc E} O(h) + O(h^2),\\
&\leq \| \widetilde{u} -u_h \|^2_{\mc E} +  O(h^{1+p}), \mbox{ using \eqref{eq:2}},\\
&\leq \|\widetilde{u} \|_{\mc {\widetilde{E}}}^2 - \| u_h \|_{\mc {\widetilde{E}}}^2 +  O(h^{1+p}), \mbox{ using \eqref{eq:4}}.
\end{align*}

Finally, using the first part of the proof: $\left| \|u\|^2_{\mc E} - \|\widetilde{u}\|^2_{\mc {\widetilde{E}}}\right| = O(h^2)$, it holds by triangular inequality
\begin{align*}
\| u -u_h \|^2_{\mc E} &\leq \left|\|\widetilde{u} \|_{\mc {\widetilde{E}}}^2 - \|u\|^2_{\mc E}\right| + \left| \|u\|^2_{\mc E}- \| u_h \|_{\mc {\widetilde{E}}}^2\right| +  O(h^{1+p}),\\
&\leq \left|\|u\|^2_{\mc E} - \| u_h \|_{\mc {\widetilde{E}}}^2\right| +  O(h^{1+p}),
\end{align*}
which shows the first inequality since $1+p = 3/2$ for FEM and 2 for GFEM \& SGFEM.

For the second inequality, we begin with
\begin{align*}
\left|\|u\|^2_\mc{E}-\|u_h\|^2_\mc{\widetilde{E}}\right| &= \left|\|u\|^2_\mc{E} -\|\widetilde{u}\|^2_{\mc {\widetilde{E}}} +\|\widetilde{u}\|^2_{\mc {\widetilde{E}}} -\|u_h\|^2_\mc{\widetilde{E}}\right|,\\
&\leq \left|\|u\|^2_\mc{E} -\|\widetilde{u}\|^2_{\mc {\widetilde{E}}}\right| +\|\widetilde{u}\|^2_{\mc {\widetilde{E}}} -\|u_h\|^2_\mc{\widetilde{E}},\mbox{ by triangular inequality},\\
&\leq \|\widetilde{u}-u_h\|^2_\mc{\widetilde{E}} + O(h^2),
\end{align*}
where again, we have used the first part of the proof: $\left| \|u\|^2_{\mc E} - \|\widetilde{u}\|^2_{\mc {\widetilde{E}}}\right| = O(h^2)$ and then Galerkin orthogonality between $\widetilde{u}-u_h$ and $u_h$ in the $\widetilde{B}(\cdot,\cdot)$ inner-product. Now, using the same kind of reasoning as for the first inequality, we have
\begin{align*}
\left|\|u\|^2_\mc{E}-\|u_h\|^2_\mc{\widetilde{E}}\right| &\leq \|\widetilde{u}-u_h\|^2_\mc{\widetilde{E}} + O(h^2),\\
&\leq \|\widetilde{u}-u_h\|^2_\mc{E} + O(h^2), \mbox{ using \eqref{eq:3}},\\
&\leq \left(\|u-\widetilde{u}\|_\mc{E} + \|u-u_h\|_\mc{E}\right)^2 + O(h^2), \mbox{ using the triangular inequality},\\
&\leq \|u-\widetilde{u}\|^2_\mc{E} + \|u-u_h\|^2_\mc{E} + 2\|u-\widetilde{u}\|_\mc{E} \|u-u_h\|_\mc{E} + O(h^2),\\
&\leq \|u-u_h\|^2_\mc{E} + O(h) \|u-u_h\|_\mc{E} + O(h^2),\mbox{ using \eqref{eq:equiv} with $p=2$}.
\end{align*}
Now, by triangular inequality, it holds
\begin{align*}
\|u-u_h\|_\mc{E} &\leq \|u-\widetilde{u}\|_\mc{E} + \|\widetilde{u}-u_h\|_\mc{E},\\
&\leq \|\widetilde{u}-u_h\|_\mc{E} + O(h),\mbox{ using \eqref{eq:equiv} with $p=2$}.
\end{align*}
So that
\begin{align*}
\left|\|u\|^2_\mc{E}-\|u_h\|^2_\mc{\widetilde{E}}\right| &\leq  \|u-u_h\|^2_\mc{E} + O(h) \|\widetilde{u}-u_h\|_\mc{E} + O(h^2).
\end{align*}
Using now \eqref{eq:2}, it yields
\begin{align*}
\left|\|u\|^2_\mc{E}-\|u_h\|^2_\mc{\widetilde{E}} \right|&\leq  \|u-u_h\|^2_\mc{E} +O(h^{1+p}),
\end{align*}
which ends the proof, since $1+p = 3/2$ for FEM and 2 for GFEM \& SGFEM.

\section{Derivation of the stopping criteria for the iterative solvers}
\label{app:stoppingcrit}
In this last annex, we derive the stopping criteria for the iterative solvers discussed in Section \ref{sec6}. The stopping criterion shown in \eqref{eq:estimfem} is derived as follows.
First, let us assume, as indicated by a priori error estimation, that there exists a constant $A>0$, independent of the mesh size $h$, such that
\begin{align*}
h^p \leq A \epsilon^h,
\end{align*}
where $\epsilon^h$ is the discreti«ation error, $p = 1/2$ for FEM and $p = 1$ for GFEM \& SGFEM.

 There exist \cite{Johnson2012} constants $B,C$, independent of the mesh size $h$, such that the following inverse estimates hold for all $v_h = \sum_{k\in \mc{N}^h_d} c_kN_k \in S^h$, denoting $\bm{c}  =[c_k]_{k\in \mc{N}^h_d}$
\begin{align*}
\| v_h \|_{L^2(\Omega)} &\leq B h \| \bm{c} \|_{l^2},\\
\| v_h \|_{\mc{E}} &\leq C h^{-1} \| v_h \|_{L^2(\Omega)}.
\end{align*}
So it holds
\begin{align*}
\| v_h \|_{\mc{E}} &\leq BC \| \bm{c} \|_{l^2}.
\end{align*}
Now, at some iteration $i$, we have an approximate solution $v_h^i =  \sum_{k\in \mc{N}^h_d} c_k^iN_k$ to \eqref{goodLinSys}, and we can form the residual vector $\bm{r}^i = \bm{f} -\bm{A}_{11}\bm{c}^i$. Of course, the exact (discrete) solution $u_h = \sum_{k\in \mc{N}^h_d} c_kN_k$ solves the (discrete) residual equation and we have $\bm{r}^i = \bm{A}_{11}(\bm{c} - \bm{c}^i)$. Using the spectral radius of $\bm{A}_{11}^{-1}$, we have
\begin{align*}
\| \bm{c} - \bm{c}^i \|_{l^2} &= \| \bm{A}_{11}^{-1} \bm{r}^i \|_{l^2},\\
&\leq \rho(\bm{A}_{11}^{-1}) \|\bm{r}^i \|_{l^2}.
\end{align*}
Further, since the condition number of $\bm{A}_{11}$ follows $\kappa_2 (\bm{A}_{11})  = O(h^{-2})$ (we mention that the largest eigenvalue of $\bm{A}_{11}$ is bounded independently of $h$), there exists a constant $D$, independent of the mesh size $h$, such that
\begin{align*}
\rho(\bm{A}_{11}^{-1}) \leq \frac{D}{h^2}.
\end{align*}
As a result, if we compute $e^i =\displaystyle \frac{\|\bm{f} -\bm{A}_{11}\bm{c}^i\|_{l^2}}{h^2}$ and perform iterations until $e^i < \epsilon/k$, where $\epsilon = h^{1/2}$ (recall that this is the FEM case), we obtain
\begin{align*}
\delta_i &= \| u_h - v_h^i\|_{\mc{E}},\\
&\leq BC \| \bm{c} - \bm{c}^i \|_{l^2},\\
&\leq BC \rho(\bm{A}_{11}^{-1}) \|\bm{r}^i \|_{l^2},\\
&\leq \frac{BCD}{h^2} \|\bm{f} -\bm{A}\bm{c}^i \|_{l^2},\\
&\leq BCD e^i,\\
&< \frac{BCD}{k}\epsilon,\\
&< \frac{BCD}{k} h^{1/2},\\
&< \frac{ABCD}{k} \epsilon^h.
\end{align*}
Thus, iterations are performed until the truncation error $\delta_i$ is smaller than the discretization error $\epsilon^h$, up to a proportionality factor controlled by $k$ and independent of the mesh size $h$. The constants $A,B,C$ and $D$ are unknown and thus taking $k$ large enough, $\delta_i$ could be made sufficiently smaller than the discretization error.

We further note that in practice, the ``effective" condition number of $\bm A_{11}$ is reduced from $O(h^{-2})$ to $O(h^{-1})$ thanks to the FMG preconditioner. As a result, it is sufficient to compute $e^i =\displaystyle \frac{\|\bm{f} -\bm{A}_{11}\bm{c}^i\|_{l^2}}{h}$ and perform iterations until $e^i < \epsilon/k$, where $\epsilon = h^{1/2}$.

The stopping criterion shown in \eqref{eq:gscheat} is derived as follows. First, since the outside iteration scheme yields a geometrical decrease of the truncation error $\delta_i = \|u_h - v_h^i \|_{\mc{E}}$, we make the following assumption: there exist positive constants $B_2 \geq B_1  \geq 0$ and  $0 <  q < 1$, all independent of $i$ such that
\begin{align*}
B_1 q^i \leq \delta_i \leq B_2 q^i.
\end{align*}
We further assume that these bounds are somewhat sharp and thus do not overlap from one iteration to the next, that is
\begin{align}
B_2 q < B_1. \label{sharpbound}
\end{align}
This is only required so that the truncation error effectively decreases at each iteration: $\delta_{i+1} < \delta_i$, which is the behavior observed in our numerical experiments.

Then, it follows by triangular inequality
\begin{align*}
\|v_h^i - v_h^{i-1} \|_{\mc{E}} &\leq \|u_h - v_h^i \|_{\mc{E}} + \|u_h - v_h^{i-1} \|_{\mc{E}},\\
&\leq B_2 (1+q) q^{i-1},
\end{align*}
and
\begin{align*}
\|v_h^i - v_h^{i-1} \|_{\mc{E}} &\geq  \| u_h- v_h^{i-1} \|_{\mc{E}} - \|u_h-v_h^i\|_{\mc{E}},\\
&\geq (B_1 - q B_2) q^{i-1}.
\end{align*}
Thus we have bounded $\|v_h^i - v_h^{i-1} \|_{\mc{E}}$ by above and below
\begin{align*}
C_1 q^{i-1} \leq \|v_h^i - v_h^{i-1} \|_{\mc{E}} \leq C_2 q^{i-1}, 
\end{align*}
with $C_1 = B_1 - q B_2 >0$ and $C_2 = B_2 (1+q)>0$, also independent of $i$. The computable quantity $\|v_h^i - v_h^{i-1} \|_{\mc{E}}$ thus follows the same behavior with the number of iterations as the truncation error $\delta_i$. Using the last three iterates we can estimate the common ratio and the scale factor of this geometric sequence. Similarly to \eqref{sharpbound}, we will assume that there exist bounds $D_2 \geq D_1  \geq 0$ in \eqref{bound2} that are sharp and do not overlap from one iteration to the next, that is
\begin{align}
D_1 q^{i-1} \leq \|v_h^i &- v_h^{i-1} \|_{\mc{E}} \leq D_2 q^{i-1}, \mbox{ and} \label{bound2},\\
D_2 q &< D_1. \nonumber
\end{align}
Again, the second condition is only required so that at each iteration: $\|v_h^i - v_h^{i-1} \|_{\mc{E}} < \|v_h^{i-1} - v_h^{i-2} \|_{\mc{E}}$, which is also the behavior observed in our numerical experiments.

Let us now recall the stopping criterion $e^i$ of \eqref{eq:gscheat}
\begin{align*}
e^i &= \frac{1}{\displaystyle\frac{1}{\|v^i_h - v^{i-1}_h\|_{\mc{E}}}-\frac{1}{\|v^{i-1}_h - v^{i-2}_h\|_{\mc{E}}}}.
\end{align*}
Now using the bounds in \eqref{bound2}, it yields
\begin{align*}
\frac{1}{\displaystyle\frac{1}{D_1q} - \frac{1}{D_2}}q^{i-2} \leq e^i \leq \frac{1}{\displaystyle\frac{1}{D_2q} - \frac{1}{D_1}}q^{i-2}.
\end{align*}
Iterations are performed until $e^i < \epsilon/k$, where $\epsilon = h$ (recall that these are the GFEM \& SGFEM cases). When this happens, we have
\begin{align*}
\delta_i &\leq B_2 q^i,\\
&\leq B_2 \left( \frac{1}{D_1q} - \frac{1}{D_2} \right)q^2 e^i,\\
&< \frac{B_2q^2}{k} \left( \frac{1}{D_1q} - \frac{1}{D_2} \right) \epsilon,\\
&< \frac{B_2q^2}{k} \left( \frac{1}{D_1q} - \frac{1}{D_2} \right) h,\\
&<  \frac{AB_2q^2}{k} \left( \frac{1}{D_1q} - \frac{1}{D_2} \right) \epsilon^h.
\end{align*}
Thus, iterations are performed until the truncation error $\delta_i$ is smaller than the discretization error $\epsilon^h$, up to a proportionality factor controlled by $k$ and independent of the mesh size $h$. Again, the constants $A,B_2,D_1,D_2$ and $q$ are unknown, and thus taking $k$ large enough, $\delta_i$ could be made sufficiently smaller than the discretization error.

\bibliography{biblio}
\bibliographystyle{plain}

\end{document}